\documentclass[12pt]{article}
\usepackage{amssymb}

\newcommand{\be}{\begin{equation}}
\newcommand{\ee}{\end{equation}}

\newcommand{\bea}{\begin{eqnarray}}
\newcommand{\eea}{\end{eqnarray}}

\newcommand{\eqq}{\!\! & = & \!\!}

\newcommand{\hst}[1]{\rule{#1}{0mm}}
\newcommand{\ty}{\hspace{0.04em}}
\newcommand{\tty}{\hspace{0.06em}}

\newcommand{\sect}[1]{\setcounter{equation}{0}\section{#1}}

\newcommand{\reff}[1]{(\ref{#1})}

\newcommand{\dis}{\displaystyle}

\newcommand{\Ga}{\Gamma}
\newcommand{\De}{\Delta}

\newcommand{\al}{\alpha}
\newcommand{\bet}{\beta}
\newcommand{\ga}{\gamma}
\newcommand{\de}{\delta}

\newcommand{\ve}{\varepsilon}

\newcommand{\lam}{\lambda}
\newcommand{\vr}{\varrho}
\newcommand{\si}{\sigma}
\newcommand{\vp}{\varphi}
\newcommand{\om}{\omega}

\newcommand{\la}{\langle}
\newcommand{\ra}{\rangle}
\newcommand{\rar}{\rightarrow}
\newcommand{\lra}{\longrightarrow}
\newcommand{\ti}{\times}
\newcommand{\op}{\oplus}

\newcommand{\ot}{\otimes}

\newcommand{\KK}{\mathbb K}

\newcommand{\ZZ}{\mathbb Z}

\newcommand{\CC}{\mathbb C}
\newcommand{\II}{\mathbb I}

\newcommand{\cd}{\cdot}

\newcommand{\comp}{\!\stackrel{\textstyle{\hst{0ex} \atop \circ}}{\hst{0ex}}\!}
\newcommand{\ol}{\overline}

\newcommand{\Uq}{U_q({\mathfrak{spo}}(2n \ty|\ty 2m))}
\newcommand{\Uqtt}{U_q({\mathfrak{spo}}(2 \ty|\ty 2))}
\newcommand{\id}{\mbox{id}}

\newcommand{\Rh}{\hat{R}}
\newcommand{\Cqi}{((C^q)^{-1})}
\newcommand{\Lgr}{\mbox{Lgr}}
\newcommand{\otb}{\:\overline{\otimes}\:}
\newcommand{\SP}{\mbox{SPO}_q(2n \ty|\ty 2m)}
\newcommand{\SPu}{\mbox{SPO}_1(2n \ty|\ty 2m)}
\newcommand{\SPSO}{\mbox{SPSO}_q(2n \ty|\ty 2m)}
\newcommand{\SPSOu}{\mbox{SPSO}_1(2n \ty|\ty 2m)}

\newcommand{\bt}{\tilde{b}}

\newcommand{\tti}{\tilde{t}}
\newcommand{\tb}{\overline{t}}
\newcommand{\Vf}{\overline{V}_{\!\!4}}

\newcommand{\Hc}{H^{\circ}}
\newcommand{\Vag}{V^{\ast\rm gr}}
\newcommand{\Wm}{W^{\rm mod}}
\newcommand{\Arq}{A^r_q(n \ty|\ty m\tty; c)}

\newcommand{\ru}{{\displaystyle{\,{}^r}}}
\newcommand{\su}{{\displaystyle{\,{}^s}}}
\newcommand{\rsu}{{\displaystyle{\,{}^{rs}}}}
\newcommand{\vso}{\vspace{-1.0ex}}

\begin{document}

\begin{titlepage}

\newlength{\ppn}
\newcommand{\defboxn}[1]{\settowidth{\ppn}{#1}}
\defboxn{BONN--TH--2000--04}

%\mbox{}
\hspace*{\fill} \parbox{\ppn}{BONN--TH--2000--04 \\
                              April 2000        }

%\vspace{10mm}
\vspace{18mm}

\begin{center}
{\LARGE\bf The quantum supergroup $\SP$  \\[0.4ex]
           and an $\SP$--covariant \\[0.7ex]
           quantum Weyl superalgebra} \\
\vspace{15mm}
{\large M. Scheunert} \\
\vspace{1mm}
Physikalisches Institut der Universit\"{a}t Bonn \\
Nu{\ss}allee 12, D--53115 Bonn, Germany \\
\end{center}

\vspace{15mm}
\begin{abstract}
\noindent
Recently, the $R$--matrix of the symplecto--orthogonal quantum superalgebra
$\Uq$ in the vector representation has been calculated. In the present work,
this $R$--matrix is used to introduce the corresponding quantum supergroup
$\SP$ and to construct an $\SP$--covariant quantum Weyl superalgebra.
\end{abstract}

\vspace{\fill}
\noindent
math.QA/0004033

\end{titlepage}

\setcounter{page}{2}
\renewcommand{\baselinestretch}{1}
\small\normalsize

\sect{Introduction \vspace{-1ex}}
Once the quantum algebras and quantum groups had been introduced \cite{Dri},
\cite{Jim}, \cite{RTF}, it was obvious to apply the same techniques also in
the super case. Indeed, at present there exist a considerable number of
papers on quantum supergroups. Most of them consider deformations of the
general linear supergroup ${\rm GL}(m \ty|\ty n)$ \cite{Man}, \cite{ZGL},
but more recently there also appeared a paper by Lee and Zhang \cite{LZh},
in which a deformation of ${\rm OSP}(1|\ty 2n)$ is discussed, and a paper by
Zhang \cite{ZOS}, dealing with a deformation of ${\rm OSP}(2 \ty|\ty 2n)$.

Quantum supergroups and quantum superalgebras are dual objects. This implies
that there are (at least) two methods to introduce a quantum supergroup. One
may regard the quantum supergroup as the primary object and define it by
means of the $R$--matrix approach, thus generalizing the methods of
Ref.~\cite{RTF} to the super case. On the other hand, one may also start
from the quantum superalgebra and construct the (algebra of functions on
the) quantum supergroup as a sub--Hopf--superalgebra of the finite dual of
the former. In Manin's paper \cite{Man} the first approach is used, while
Zhang \cite{ZGL}, \cite{LZh}, \cite{ZOS} prefers to use the second. It is
believed that both methods should lead to the same objects, however, since
it is notoriously difficult to put the duality of quantum (super)groups and
quantum (super)algebras on a firm basis (including a proof that the dual
pairing is non--degenerate(!)), there is still a lot to be done to
substantiate this belief (for example, see Ref.~\cite{DTa}; this paper
contains an extensive list of references on quantum supergroups and quantum
superalgebras).

The present paper is a direct sequel to Ref.~\cite{Srm}, in which I have
calculated the $R$--matrix of the quantum superalgebra $\Uq$ in the vector
representation. Here, I am going to use this $R$--matrix to introduce the
corresponding quantum supergroup $\SP$ and to construct some of its comodule
superalgebras, among which is an $\SP$--covariant quantum Weyl superalgebra.
This will be done by means of techniques which (in the non--super case) have
been described in Ref.~\cite{RTF} (i.e., we use the first of the two
approaches described above). However, as we are going to see, in the super
case some amendments are necessary in order to cope with the nototious
problem that the representations are not necessarily completely reducible.
The reason for considering the $\SP$ quantum supergroup rather than
${\rm OSP}_q(2m \ty|\ty 2n)$ has been explained in Ref.~\cite{Srm}.

The present work falls into two parts. The first half (Sections 2 to 4)
contains results which are applicable to arbitrary quantum supergroups (not
just the symplecto--orthogonal ones). In fact, these results are formulated
in the general graded framework, where $\Ga$ in an arbitrary abelian group
and $\si$ is a commutation factor on $\Ga$. In the second half (Sections 5
and 6), the special case of the $\SP$ quantum supergroups is considered.

More precisely, this paper is set up as follows. In Section 2 we are going to
recall some elementary facts about matrix elements of graded representations
of a $\si$--bialgebra or $\si$--Hopf algebra. These facts will be used to
motivate the definitions and constructions of the subsequent sections. In
particular, they serve to find the appropriate sign factors (i.e.,
commutation factors) by which the formulae of the present paper differ from
those in the non--graded setting.

Usually, a quantum (super)group is constructed in two steps. First, one
introduces the corresponding bi--(super)algebra $A(R)$, then one amends this
definition to obtain the quantum (super)group itself, which is a Hopf
(super)algebra. Basic to the construction of $A(R)$ is a finite--dimensional
(graded) vector space $V$ and some (even) linear operator $R$ in $V \ot V$
which, in the more advanced part of the theory, is assumed to be invertible
and to satisfy the (graded) Yang--Baxter equation. The $\si$--bialgebras
$A(R)$ are discussed to some extent in Section 3. In particular, we prove
(under the foregoing additional assumptions on $R\ty$) that $A(R)$ is
coquasitriangular.

Section 4 is devoted to the discussion of graded comodules and comodule
algebras over $A(R)$. We concentrate on those of these structures that are
derived from the vector comodule or its dual.

In Section 5 we introduce the quantum supergroup $\SP$ and derive some of its
properties. In particular, we show that $\SP$ is coquasitriangular, and that
there exists a	 natural dual pairing of the Hopf superalgebras $\SP$ and
$\Uq$. (According to the general folklore, both of these results are expected
to be true.)

Section 6 contains a discussion of some graded $\SP$--comodule superalgebras.
After describing the general setting, we concentrate on the construction of
the $\SP$--covariant quantum Weyl superalgebras.

The final Section 7 contains a discussion of our results and some hints to
further research.

We close this introduction by some comments on our conventions and notation.
Throughout the present work, the base field will be an arbitrary field $\KK$
of characteristic zero. The pairing of a vector space and its dual will be
denoted by pointed brackets, and the action of an algebra $A$ on an
$A$--module will be denoted by a dot. It should be obvious that we are going
to work in the general setting of graded algebraic structures (see
Ref.~\cite{Sgt}). To remind the reader of this fact, the graded tensor
product of linear maps and of graded algebras will be denoted by $\otb$.
However, it might be worthwhile to stress that even if the overbar seems to
be missing this does {\em not}\/ indicate that we have returned to the
non--graded setting.

\sect{Some elementary properties of matrix elements of representations
                                                               \vspace{-1ex}}
In the present section we recall some elementary properties of the matrix
elements of a graded representation of a $\si$--Hopf algebra (or
$\si$--bialgebra). These properties will serve to motivate the definitions
and constructions given later.

First of all, we introduce some notation. Throughout the present section,
$\Ga$ will be an abelian group, and $\si$ will be a commutation factor on
$\Ga$ with values in $\KK\,$. All gradations appearing in this section are
understood to be $\Ga$--gradations. In the following, we shall freely use the
notation and results of Ref.~\cite{Sgt}.

Let $V$ be a {\em finite--dimensional} graded vector space, let $\Vag$ be the
graded dual of $V$, and let $\Lgr(V)$ be the space of all linear maps of $V$
into itself that can be written as a sum of homogeneous linear maps. In
Ref.~\cite{Sgt}, we have denoted the latter space by $\Lgr(V,V)$. Since $V$
is finite--dimensional, it is known that, regarded as vector spaces, $\Vag$
is equal to $V^{\ast}$, the dual of $V$, and $\Lgr(V)$ is equal to the space
of all linear maps of $V$ into itself. We use the present notation in order
to stress that $\Vag$ is a graded vector space, and that $\Lgr(V)$ is a graded
associative algebra.

Suppose now that $(e_i)_{i \in I}$ is a homogeneous basis of $V$, where $I$
is some index set, and that $e_i$ is homogeneous of degree $\eta_i \in \Ga$.
Let $(e'_i)_{i \in I}$ be the basis of $\Vag$ that is dual to
$(e_i)_{i \in I}\tty$, and let $(E_{ij})_{i,\ty j \in I}$ be the basis of
$\Lgr(V)$ that canonically corresponds to $(e_i)_{i \in I}\,$:
 \[ E_{ij}(e_k) = \de_{jk} \ty\ty e_i
                                    \quad\mbox{for all $i,j,k \in I$} \,.\] 
Then $e'_i$ is homogeneous of degree $-\eta_i\ty\ty$, and $E_{ij}$ is
homogeneous of degree $\eta_i - \eta_j\ty\ty$.

It is well-known that the class of $\Ga$--graded vector spaces, endowed with
the usual tensor product of graded vector spaces and with the twist maps
defined by means of $\si\ty$, forms a tensor category (see Ref.~\cite{Sgt} for
details). A (generalized) Hopf algebra $H$ living in this category is called
a $\si$--Hopf algebra. More explicitly, $H$ is a $\Ga$--graded associative
algebra with a unit element, and it is endowed with a coproduct $\De\ty\ty$,
a counit $\ve\ty\ty$, and an antipode $S\ty\ty$, which satisfy the obvious
axioms (in the category). In particular, this implies that the structure
maps $\De\ty\ty$, $\ve\ty\ty$, and $S$ are homogeneous of degree zero.
Needless to say, $\si$--bialgebras are introduced similarly. (For generalized
Hopf algebras living in more general categories, see Ref.~\cite{Maj}.)

Let $H$ be a $\si$--Hopf algebra, and let $\Hc$ be the set of all linear
forms on $H$ that vanish on a graded two--sided ideal of $H$ of finite
codimension. This definition generalizes the definition given in Sweedler's
book \cite{Swe} to the general graded case considered here.

It is easy to see that the properties derived in Ref.~\cite{Swe} can
immediately be extended to the present setting. In particular, this is true
for the criteria characterizing the elements of $\Hc$. We don't want to go
into detail here, but mention some properties not given there.

First of all we note that in the preceding definition the graded two--sided
ideals can be replaced by graded left or right ideals (indeed, every graded
left or right ideal of finite codimension contains a graded two--sided ideal
of finite codimension). Secondly, the elements of $\Hc$ can also be
characterized by the following property: \\[1.0ex]
A linear form $f$ on $H$ belongs to $\Hc$ if and only if there exists a
finite--dimensional graded (left) $H$--module $V$, an element $x \in V$, and
a linear form $x' \in \Vag$ such that
 \[ f(h) = \la x', h \cd x \ra \quad\mbox{for all $h \in H$} \,. \]
This property shows very clearly that $\Hc$ is closely related to the theory
of finite--dimensional graded representations of $H$.

Obviously, $\Hc$ is a graded subspace of $H^{\ast\rm gr}$, the graded dual
of $H$. Actually, much more is true, for the structure of $H$ as a
$\si$--Hopf algebra leads to a similar structure on $\Hc$. More precisely,
the transpose of $\De$ induces the multiplication in $\Hc$, the counit
$\ve$ is the unit element of $\Hc$, the transpose of the product map of
$H$ induces the coproduct $\De_{\circ}$ of $\Hc$, the evaluation at the unit
element of $H$ is the counit $\ve_{\circ}$ of $\Hc$, and the transpose of $S$
induces the antipode $S_{\circ}$ of $\Hc$. This $\si$--Hopf algebra $\Hc$ is
called the finite (or continuous) dual of $H$. 

Using the notation introduced at the beginning of this section, let us now
assume in addition that $V$ is a graded $H$--module, and that
$\pi : H \rar \Lgr(V)$ is the graded representation afforded by it. If
$i,j \in I$, we define the linear form $\pi_{ij}$ on $H$ by
 \[ \pi_{ij}(h) = \la e'_i \ty, \pi(h)\ty e_j \ra
                                       \quad\mbox{for all $h \in H$} \,. \]
Obviously, $\pi_{ij}$ is an element of $\Hc$, and is homogeneous of degree
$\eta_j - \eta_i\ty\ty$. Moreover, we have
 \[ \pi(h) = \sum_{i,\ty j \in I} \pi_{ij}(h) \ty E_{ij}
                                       \quad\mbox{for all $h \in H$} \,. \]
Then our first important observation is that
 \be \De_{\circ}(\pi_{ij})
       = \sum_{k \in I} \si(\eta_k - \eta_i \ty, \eta_k - \eta_j)
                                \tty \pi_{ik} \ot \pi_{kj} \label{copr} \ee
for all $i,j \in I$. By definition of the antipode, this implies that
 \be \sum_{k \in I} \si(\eta_k - \eta_i \ty, \eta_k - \eta_j)\tty
                                         S_{\circ}(\pi_{ik}) \tty \pi_{kj}
  \,=\, \sum_{k \in I} \si(\eta_k - \eta_i \ty, \eta_k - \eta_j)\ty\ty
                                         \pi_{ik} \ty\ty S_{\circ}(\pi_{kj})
  \,=\, \de_{i,\ty j} \ty\ty 1_{\Hc} \;, \label{Spipi} \ee
for all $i,j \in I$, where we have used the obvious fact that
 \[ \ve_{\circ}(\pi_{ij}) = \de_{ij} \quad\mbox{for all $i,j \in I$} \,.  \]

Let us next recall that $V$ has a canonical structure of a graded right
$\Hc$--comodule. It is uniquely fixed by the requirement that it induces, in
a canonical way, the original structure of a graded left $H$--module on $V$.
More precisely, let 
 \[ \de : V \lra V \ot \Hc \]
be the structure map of $V$, regarded as a graded right $\Hc$--comodule, let
$x \in V$ be homogeneous of degree $\xi\tty$, and set
 \[ \de(x) = \sum_\ru x^0_r \ot x^1_r \;, \vso \]
with $x^0_r \in V$ and $x^1_r \in \Hc$. Then, if $h \in H$ is homogeneous of
degree $\eta\tty$, we have
 \[ h \cd x = \pi(h)\tty x = \si(\eta,\xi) \sum_\ru
                                   x^0_r \ty\la x^1_r \ty, h \ra \;, \vso \]
where $\la \;,\: \ra$ denotes the dual pairing of $\Hc$ and $H$. Explicitly,
we find
 \be \de(e_i) = \sum_{j \in I} \si(\eta_i \tty, \eta_j - \eta_i) \tty
                                       e_j \ot \pi_{ji} \;, \label{rcom} \ee
for all $i \in I$.

Similarly, let $V'$ be the graded {\em right} $H$--module that is dual to
$V$. Regarded as a graded vector space, $V'$ is equal to $\Vag$, and the
module structure is given by
 \[ \la x' \cd h \ty\ty, x \ra = \la x', h \cd x \ra \;, \]
for all $h \in H$, $x \in V$, and $x' \in V'$. (Of course, this time
$\la \;,\: \ra$ denotes the dual pairing of $\Vag$ and $V$.) More explicitly,
we have
  \[ e'_i \cd h = \sum_{j \in I} \pi_{ij}(h)\tty e'_j \;, \]
for all $h \in H$ and all $i \in I$. Then $V'$ has a canonical structure of
a graded left $\Hc$--comodule. It is uniquely fixed by the requirement that
it induces, in a canonical way, the original structure of a graded right
$H$--module on $V'$. More precisely, let 
 \[ \de' : V' \lra \Hc \ot V' \]
be the structure map of $V'$, regarded as a graded left $\Hc$--comodule, let
$z \in V'$, and set
 \[ \de'(z) = \sum_\ru z^{-1}_r \ot z^0_r \;, \vso \]
where $z^{-1}_r \in \Hc$, and where $z^0_r \in V'$ is homogeneous of degree
$\zeta^0_{\ty r}\tty$. Then, if $h \in H$ is homogeneous of degree
$\eta\tty$, we have
  \[ z \cd h = \sum_\ru \si(\zeta^0_{\ty r} \ty, \eta) \ty
                                \la z^{-1}_r , h \ra \tty z^0_r \;. \vso \]
Explicitly, we find
  \be \de'(e'_i) = \sum_{j \in I} \si(\eta_j \tty, \eta_i - \eta_j) \tty
                                     \pi_{ij} \ot e'_j \;, \label{lcom} \ee
for all $i \in I$.

Let us next consider two cases where our assumptions are more restrictive.
In the first case, we suppose that on $V$ there exists a non--degenerate
invariant bilinear form $b$ which is homogeneous of degree zero (see
Appendix A of Ref.~\cite{Srm} for some comments on invariant bilinear forms).
Define a linear map
 \[ f_r : V \lra \Vag \]
by
  \[ \la f_r(y), x \ra = \si(\eta\ty,\xi)\ty b(x,y) \;, \]
where $x$ and $y$ are elements of $V$ which are homogeneous of degrees $\xi$
and $\eta\tty$, respectively. Then the assumption that $b$ is invariant is
equivalent to the requirement that
 \[ \pi(S(h)) = f_r^{-1} \comp (\pi(h))^{\rm T} \comp f_r
                                       \quad\mbox{for all $h \in H$} \,, \]
where ${}^{\rm T}$ denotes the $\si$--transpose of a linear map (see
Eqn.~(A.5) of Ref.~\cite{Srm}). In terms of matrix elements, this condition
takes the form
 \be \pi_{ij}(S(h)) = \sum_{k,\ty\ell} \ty
                                (F_r^{-1})_{ik} \,
                                     (\pi(h)^{\rm T})_{k\ell} \,\ty
                                                  (F_r)_{\ell j}
    \quad\mbox{for all $i,j \in I$ and all $h \in H$} \,, \label{piS} \ee
where $F_r$ is the matrix of $f_r\,$:
 \[ f_r(e_j) = \sum_{i \in I} (F_r)_{ij} \tty\ty e'_i
                                       \quad\mbox{for all $j \in I$} \,. \]
To understand the meaning of Eqn.~\reff{piS} we note that
 \be \pi_{ij} \comp S = S_{\circ}(\pi_{ij}) \;, \label{Scpi} \ee
moreover, we recall that the matrix elements of $\pi(h)^{\rm T}$ are given
by
 \be (\pi(h)^{\rm T}\ty)_{k\ell}
              = \si(\eta_{\ell}\ty, \eta_{\ell} - \eta_k)
                                 \ty\ty \pi_{\ell k}(h) \;. \label{piT} \ee
Thus Eqn.~\reff{piS} expresses $S_{\circ}(\pi_{ij})$ as a linear combination
of the matrix elements $\pi_{\ell k}\ty\ty$.

For convenience, we also give a formula for the matrix $F_r\ty\ty$. Setting
 \[ C_{ij} = b(e_i \ty, e_j) \]
and
 \[ \si_{ij} = \si(\eta_i \ty, \eta_j) \quad,\quad
                                    \si_i =  \si(\eta_i \ty, \eta_i) \;, \]
for all $i,j \in I$, we obtain
 \be (F_r)_{ij} = \si_{ji} \ty\ty C_{ij}
                   = \si_i \ty\ty C_{ij}
                      = C_{ij} \tty\ty \si_j 
                         \quad\mbox{for all $i,j \in I$} \,. \label{Fr} \ee

In the second case, we assume in addition that $H$ is quasitriangular, i.e.,
that we are given a universal $R$--matrix $\cal R$ for $H$. By definition,
$\cal R$ is an invertible element of $H \otb H$ and is homogeneous of degree
zero. Moreover, it satisfies the following relations:
 \be {\cal R} \tty \De(h) = (P \De(h)) \tty {\cal R} \quad
               \mbox{for all $h \in H$} \,, \vspace{0.5ex} \label{Rdel} \ee
 \be (\De \ot \id_H)({\cal R}) = {\cal R}_{13} {\cal R}_{23} \quad,\quad
     (\id_H \ot \De)({\cal R}) = {\cal R}_{13} {\cal R}_{12}
                                        \vspace{0.5ex} \label{hexa} \;. \ee
Our notation is standard: $P : H \otb H \lra H \otb H$ is the twist map (in
the graded sense), and the elements ${\cal R}_{12}\tty$, ${\cal R}_{13}\tty$, 
and ${\cal R}_{23}\tty$ of $H \otb H \otb H$ are defined by
 \[ {\cal R}_{12} = {\cal R} \ot 1 \quad,\quad 
                                           {\cal R}_{23} = 1 \ot {\cal R} \]
 \[ {\cal R}_{13}
        = (P \ot \id_H)({\cal R}_{23}) = (\id_H \ot P)({\cal R}_{12}) \;. \]

Let us define the $R$--matrix $R$ (more precisely, the $R$--operator $R\ty$)
acting in $V \ot V$ by
 \[ R = (\pi \ot \pi)(\cal R) \;. \]
Using the notation
 \be {\cal R} = \sum_\ru R^1_r \ot R^2_r \;, \vso \label{Rexpl} \ee
we have
 \[ R = \sum_{ijk\ty\ell \in I} \!\!
                     R_{ij,\ty k\ty\ell}\ty\ty E_{ik} \otb E_{j\ell} \;, \]
with
 \[ R_{ij,\ty k\ty\ell} = \sum_\ru \pi_{ik}(R^1_r) \, \pi_{j\ell}(R^2_r) \;.
                                                                      \vso \]
Acting with $\pi \ot \pi$ on both sides of Eqn.~\reff{Rdel}, we obtain
 \be \sum_{a,b \in I} \si(\eta_j - \eta_{\ell} \ty, \eta_a - \eta_k)
                   \tty R_{ij,\ty ab} \ty\ty \pi_{ak} \ty\ty \pi_{b\ty\ell}
   \,=\, \sum_{a,b \in I} \si(\eta_j - \eta_b \ty, \eta_a - \eta_k)
           \ty\ty \pi_{jb} \ty\ty \pi_{ia}\ty\ty R_{ab,\ty k\ty\ell} \;,
                                                          \label{RTTex} \ee
for all $i,j,k,\ell \in I$. These equations generalize the famous
RTT--{\tty}relations. By using the following technical device, they can be
written in the familiar form.

Consider the algebra
 \[ \Lgr(V) \otb \Lgr(V) \otb \Hc \simeq \Lgr(V \ot V) \otb \Hc \]
and define the following elements therein:
 \be \pi_1 \,=\, \sum_{i,\ty j \in I} E_{ij} \ot \id_V \ot \pi_{ij}
                         \quad,\quad
     \pi_2 \,=\, \sum_{i,\ty j \in I} \id_V \ot E_{ij} \ot \pi_{ij} \;.
                                                       \label{pionetwo} \ee
Then the Eqns.~\reff{RTTex} take the form
 \be (R \ot \ve) \tty \pi_1 \tty \pi_2
                     = \pi_2 \ty\ty \pi_1 \tty (R \ot \ve) \label{RTT} \ee
(recall that $\ve$ is the unit element of $\Hc\ty$). \vspace{2.0ex}

\noindent
{\em Remark 2.1.} Some comments on the notation are in order. If $A$ and
$B$ are two graded associative algebras, we denote their graded tensor
product by $A \otb B$ and the decomposable tensors in $A \otb B$ by $a \ot b$
(where $a \in A$ and $b \in B\ty$, see the appendix of Ref.~\cite{Sin}). Let
us next recall that there exists a canonical isomorphism of graded algebras
 \[ \Lgr(V) \otb \Lgr(V) \lra \Lgr(V \ot V) \;; \]
for all $g,g' \in \Lgr(V)$, it maps the tensor $g \ot g'$ onto the graded
tensor product $g \otb g'$ of the linear maps $g$ and $g'$ (see the Eqns.~(3.4)
and (3.5) of Ref.~\cite{Sgt}\ty). Originally, $R$ is an element of
$\Lgr(V \ot V)$. However, occasionally (as above) it is useful to regard it
also as an element of $\Lgr(V) \otb \Lgr(V)$. As such, it will be written in
the form
$R = \sum_{ijk\ty\ell \in I} R_{ij,\ty k\ty\ell}\ty\ty E_{ik} \ot E_{j\ell}\,$.
The reader should be careful not to confound $E_{ik} \ot E_{j\ell}$ with the
usual (non--graded) tensor product of the linear maps $E_{ik}$ and
$E_{j\ell}$ (which has been used in Ref.~\cite{Srm}, but will never be used
in the present work). \vspace{2.0ex}

Now let $A$ be a graded sub--$\si$--bialgebra of $\Hc$. Then the universal
$R$--matrix yields a bilinear form $\vr$ on $A$ that is defined by
 \[ \vr(a \ty,b) = \la a \ot b \ty, {\cal R} \ra
                                     \quad\mbox{for all $a,b \in A$} \;. \]
Explicitly, we have
 \[ \vr(a \ty,b) = \sum_\ru \si(\rho^1_r \ty, \rho^2_r)
                          \la a \ty, R^1_r \ra \la b \ty, R^2_r \ty\ra
                                \quad\mbox{for all $a,b \in A$} \;. \vso \]
(We use the notation introduced in Eqn.~\reff{Rexpl} and assume that $R^1_r$
and $R^2_r$ are homogeneous of degrees $\rho^1_r$ and $\rho^2_r\tty$,
respectively.) In particular, if the linear forms $\pi_{ij}$ belong to
$A\tty$, we have
 \be \vr(\pi_{ik} \ty, \pi_{j\ell})
     = \si(\eta_i - \eta_k \ty, \eta_j - \eta_{\ell}) R_{ij,\ty k\ty\ell}
             \quad\mbox{for all $i,j,k,\ell \in I$} \,. \label{rhopipi} \ee
Moreover, the defining properties of a universal $R$--matrix imply the
following properties of $\vr\tty$. First of all, the bilinear form $\vr$ is
homogeneous of degree zero. Secondly, $\vr$ is convolution invertible.

To understand this statement, we recall that $A \ot A$ has a canonical
structure of a graded coalgebra. Correspondingly, its graded dual
$(A \ot A)^{\ast\rm gr}$ is a graded algebra (whose multiplication is called
convolution). Since homogeneous bilinear forms on $A$ can be canonically
identified with homogeneous linear forms on $A \ot A\tty$, it makes sense
to require that $\vr$ is invertible in $(A \ot A)^{\ast\rm gr}$, and this is
the meaning of the statement above. More explicitly, this is to say that
there exists a (necessarily unique) bilinear form $\vr'$ on $A\tty$, which is
homogeneous of degree zero and satisfies
 \be \sum_{\rsu} \si(\al^2_r \ty, \bet^1_s)
                            \vr(a^1_r \ty, b^1_s) \vr'(a^2_r \ty, b^2_s) 
  \,=\, \sum_{\rsu} \si(\al^2_r \ty, \bet^1_s)
                            \vr'(a^1_r \ty, b^1_s) \vr(a^2_r \ty, b^2_s) 
  \,=\, \ve(a) \ty \ve(b)                           \vso \label{rhoinv} \ee
for all $a,b \in A\tty$. In these equations, we have set
 \[ \De(a) = \sum_\ru a^1_r \ot a^2_r \quad,\quad
                              \De(b) = \sum_\su b^1_s \ot b^2_s \;, \vso \]
and the elements $a^1_r\ty, a^2_r\ty, b^1_s\ty, b^2_s \in A$ are supposed to
be homogeneous of the degrees 
$\al^1_r\ty, \al^2_r\ty, \bet^1_s\ty, \bet^2_s\tty$, respectively.

Using the same notation, the relation \reff{Rdel} implies that
 \be    \sum_{\rsu} \si(\al^2_r \ty, \bet^1_s)
                              \vr(a^1_r \ty, b^1_s) \tty a^2_r \ty\ty b^2_s
  \,=\, \sum_{\rsu} \si(\al^1_r + \al^2_r \ty, \bet^1_s) \tty
                              b^1_s \ty\ty a^1_r \tty \vr(a^2_r \ty, b^2_s)
                                                    \vso \label{comrel} \ee
for all $a,b \in A\tty$, and the relations \reff{hexa} imply that
 \be \vr(a \tty a',b)
     = \sum_\su \si(\al', \bet^1_s) \vr(a\ty, b^1_s) \vr(a', b^2_s\ty)
                                                     \vso \label{multl} \ee
 \be \vr(a \ty, b\ty b') = \sum_\ru \vr(a^1_r\ty, b') \vr(a^2_r\ty, b)
                                                     \vso \label{multr} \ee
for all $a\ty,a',b\ty,b' \in A\tty$, where $a'$ is supposed to be homogeneous
of degree $\al'$.

A bilinear form $\vr$ on a $\si$--bialgebra $A$ that has the properties
given above is called a {\em universal $r$--form} for $A\tty$, and a
$\si$--bialgebra (or $\si$--Hopf algebra) endowed with a universal $r$--form
is said to be {\em coquasitriangular.}

If $\vr$ is a universal $r$--form for a bialgebra $A\tty$, then so is
$\vr' \comp P$, where $P : A \otb A \lra A \otb A$ denotes the twist. The
conditions above imply that
 \be \vr(a \tty, 1) = \vr(1\ty,\ty a) = \ve(a)
                       \quad\mbox{for all $a \in A$} \;. \label{rhounit} \ee
(To prove this one has to use the fact that $\vr$ is invertible.) Moreover,
if $A$ is a $\si$--Hopf algebra, it follows that
 \be \vr(S(a),b) = \vr'(a\ty,b) \quad,\quad
                   \vr'(a\ty, S(b)) = \vr(a \ty, b) \;, \label{rhoanti} \ee
for all $a,b \in A\tty$.

Let $A^{\rm cop}$ (resp.~$A^{\rm aop}$) be the $\si$--bialgebra which,
regarded as a graded algebra (resp.~as a graded coalgebra) is equal to
$A\tty$, but whose coalgebra (resp.~algebra) structure is opposite (in the
graded sense) to that of $A\tty$. Then the Eqns.~\reff{multl}, \reff{multr},
and \reff{rhounit} show that $\vr$ is a dual pairing (possibly degenerate)
of the $\si$--bialgebras $A^{\rm cop}$ and $A\tty$, and also of the
$\si$--bialgebras $A$ and $A^{\rm aop}$.

\sect{The $\si$--bialgebras $A(R)$ \vspace{-1ex}}
We keep the notation introduced at the beginning of Section 2 and start by
recalling a simple fact (see Ref.~\cite{Swe} for the non--graded case). Let
$C$ be a graded coalgebra. Equivalently, this means that $C$ is a graded
vector space and a coalgebra, and that the structure maps of the latter are
homogeneous of degree zero. Consider the tensor algebra $T(C)$ of the graded
vector space $C$. It is well--known that $T(C)$ is a $\ZZ \ti \Ga$--graded
algebra. Then there exists a unique coalgebra structure on $T(C)$ that
converts $T(C)$ into a $\si$--bialgebra, and such that the canonical
injection $\iota : C \rar T(C)$ is a homomorphism of graded coalgebras. The
$\si$--bialgebra $T(C)$ and the injection $\iota$ have the following
universal property: \\[1.0ex]
If $B$ is any $\si$--bialgebra, and if $\vp : C \rar B$ is a homomorphism of
graded coalgebras, there exists a unique $\si$--bialgebra homomorphism
$\ol{\vp} : T(C) \rar B$ such that $\vp = \ol{\vp} \comp \iota\ty\ty$.

We apply the foregoing construction to the following special case. Let the
index set $I$ and the family of degrees $\eta = (\eta_i)_{i \in I}$ be as at
the beginning of Section 2. We consider the following graded coalgebra
$C(\eta)$\,: It has a basis $(X_{ij})_{i,\ty j \in I}\tty$, such that $X_{ij}$
is homogeneous of degree $\eta_j - \eta_i\tty$, and such that the coproduct
$\De$ and the counit $\ve$ are given by
 \[ \De(X_{ij}) = \sum_{k \in I} \si(\eta_k - \eta_i \ty, \eta_k - \eta_j)
                                                   \tty X_{ik} \ot X_{kj} \]
 \[ \ve(X_{ij}) = \de_{ij} \;, \vspace{0.5ex} \]
for all $i,j \in I$. Obviously, according to Eqn.~\reff{copr}, the graded
coalgebra $C(\eta)$ is (isomorphic to) the dual of the graded algebra
$\Lgr(V)$. We stress that $T(C(\eta))$ is a free algebra with free generators
$X_{ij}\ty\ty$; $i,j \in I$.

Now let
 \[ R : V \ot V \lra V \ot V \]
be a linear map which is {\em homogeneous of degree zero.} As before, we use
the notation
 \be R = \sum_{ijk\ty\ell \in I} \!\!
         R_{ij,\ty k\ty\ell}\ty\ty E_{ik} \otb E_{j\ell} \;. \label{Rnot} \ee
The homogeneity condition is equivalent to the requirement that
 \[ R_{ij,\ty k\ty\ell} = 0
           \quad\mbox{if $\eta_i + \eta_j \neq \eta_k + \eta_{\ell}$} \;. \]
Consider the graded algebra
 \[ \Lgr(V) \otb \Lgr(V) \otb T(C(\eta))
                                     \simeq \Lgr(V \ot V) \otb T(C(\eta)) \]
and introduce the following elements therein:
 \be X_1 \,=\, \sum_{i,\ty j \in I} E_{ij} \ot \id_V \ot X_{ij}
                  \quad,\quad 
     X_2 \,=\, \sum_{i,\ty j \in I} \id_V \ot E_{ij} \ot X_{ij} \ee
(see the Eqns.~\reff{pionetwo}). Define the elements
$X_{ij,\ty k\ty\ell} \in T(C(\eta))$\,; $i,j,k,\ell \in I$, through the
following equation:
 \[ (R \ot 1) \tty X_1 \tty X_2 - X_2 \ty\ty X_1 \ty (R \ot 1)
            \,=\, \sum_{ijk\ty\ell \in I}
                             E_{ik} \ot E_{j\ell} \ot X_{ij,\ty k\ty\ell} \]
(see Eqn.~\reff{RTT}). According to the definition, $X_{ij,\ty k\ty\ell}$ is
homogeneous of degree $\eta_k - \eta_i + \eta_{\ell} - \eta_j\tty$ with
respect to the $\Ga$--gradation and homogeneous of degree 2 with respect to
the $\ZZ$--gradation. Explicitly, we obtain
 \[ X_{ij,\ty k\ty\ell}
    \,=\, \sum_{a,b \in I} \si(\eta_j - \eta_{\ell} \ty, \eta_a - \eta_k)
                   \tty R_{ij,\ty ab} \ty\ty X_{ak} \ty\ty X_{b\ty\ell}
        - \sum_{a,b \in I} \si(\eta_j - \eta_b \ty, \eta_a - \eta_k)
                \ty\ty X_{jb} \ty\ty X_{ia}\ty\ty R_{ab,\ty k\ty\ell} \;, \]
for all $i,j,k,\ell \in I$.

Let $J(R)$ be the (two--sided) ideal of $T(C(\eta))$ that is generated by the
elements $X_{ij,\ty k\ty\ell}$\,; $i,j,k,\ell \in I$. It is not difficult to
see that $J(R)$ is, in fact, a $\Ga$--graded biideal and a
$\ZZ \ti \Ga$--graded ideal of $T(C(\eta))$.  Consequently, the quotient
 \[ A(R) = T(C(\eta))/J(R) \]
inherits from $T(C(\eta))$ the structure of a $\si$--bialgebra and also
of a $\ZZ \ti \Ga$--graded algebra. (Our notation is somewhat incomplete
since, obviously, $A(R)$ also depends on $\si$ and on
$\eta = (\eta_i)_{i \in I}\tty$.) The canonical image of $X_{ij}$ in $A(R)$
will be denoted by $t_{ij}\tty$, for all $i,j \in I$.

Consider the graded algebra
 \[ \Lgr(V) \otb \Lgr(V) \otb A(R) \simeq \Lgr(V \ot V) \otb A(R) \]
and introduce the following elements therein:
 \be T_1 \,=\, \sum_{i,\ty j \in I} E_{ij} \ot \id_V \ot t_{ij}
                  \quad,\quad 
    T_2 \,=\, \sum_{i,\ty j \in I} \id_V \ot E_{ij} \ot t_{ij} \;.
                                                            \label{ToTt} \ee
Then $A(R)$ is the universal graded algebra, generated by elements
$t_{ij}\tty$; $i,j \in I$, which are homogeneous of degree $\eta_j - \eta_i$
and satisfy the following RTT--{\tty}relation
 \be (R \ot 1) T_1 T_2 = T_2 \ty T_1 (R \ot 1) \;. \label{RTTA} \ee
Explicitly, this relation is equivalent to
 \be \sum_{a,b \in I} \si(\eta_j - \eta_{\ell} \ty, \eta_a - \eta_k)
                   \tty R_{ij,\ty ab} \ty\ty t_{ak} \ty\ty t_{b\ty\ell}
   \,=\, \sum_{a,b \in I} \si(\eta_j - \eta_b \ty, \eta_a - \eta_k)
           \ty\ty t_{jb} \ty\ty t_{ia}\ty\ty R_{ab,\ty k\ty\ell} \;,
                                                          \label{RTTAex} \ee
for all $i,j,k,\ell \in I$. Moreover, the coproduct $\De$ and the counit
$\ve$ are uniquely fixed by the equations
 \be \De(t_{ij}) = \sum_{k \in I} \si(\eta_k - \eta_i \ty, \eta_k - \eta_j)
                                    \tty t_{ik} \ot t_{kj} \label{delT} \ee
 \be \ve(t_{ij}) = \de_{ij} \;,              \vspace{0.5ex} \label{epsT} \ee
for all $i,j \in I$.

It is well--known that there exists a useful equivalent formulation of
these relations. Let $P : V \ot V \rar V \ot V$ be the twist.  Explicitly, we
obtain
 \[ P = \sum_{i,\ty j \in I} \si(\eta_i \ty, \eta_i) E_{ji} \otb E_{ij}  \]
(the factor $\si(\eta_i \ty, \eta_i)$ is not a misprint). If $g$ and $g'$ are
two elements of $\Lgr(V)$ which are homogeneous of degrees $\ga$ and
$\ga\ty'$, respectively, we have
 \[ P \comp (g \otb g') = \si(\ga,\ga\ty') \ty (g' \otb g) \comp P \;. \]
Define
 \[ \Rh = P R \;. \]
Then the relations \reff{RTTA} can also be written in the form
 \be (\Rh \ot 1) T_1 T_2 = T_1 \ty T_2 (\Rh \ot 1) \;. \label{RTTbr} \ee

Up to now, $R$ could be an arbitrary element of $\Lgr(V \ot V)$ that is
homogeneous of degree zero. This indicates that we haven't yet proved
anything really substantial about the $\si$--bialgebras $A(R)$. Presumably,
the most important property of the $A(R)$'s is that, under certain natural
additional assumptions, they are coquasitriangular. More precisely, we have
the following theorem. \vspace{1.0ex}

\noindent
{\bf Theorem 1.} We keep the notation used above. Let $R$ be an element of
$\Lgr(V \ot V)$ that is homogeneous of degree zero. Suppose in addition that
$R$ is invertible and satisfies the graded Yang--Baxter equation. Then there
exists a unique universal $r$--form $\vr$ on $A(R)$ such that
 \be \vr(t_{ik} \ty, t_{j\ell})
     = \si(\eta_i - \eta_k \ty, \eta_j - \eta_{\ell}) R_{ij,\ty k\ty\ell}
                \quad\mbox{for all $i,j,k,\ell \in I$} \,. \label{rhott} \ee
(We are using the notation introduced in Eqn.~\reff{Rnot}.) \vspace{1.0ex}

\noindent
Sketch of the proof \\[1.0ex]
We start from the following fact which, in the non-graded case, has been
mentioned in Ref.~\cite{TTa}. Consider two graded coalgebras $C$ and $D$
and the corresponding $\si$--bialgebras $T(C)$ and $T(D)$ introduced at the
beginning of this section. If $\psi_1 : D \ti C \rar \KK$ is any bilinear form
which is homogeneous of degree zero, there exists a unique dual pairing $\psi$ 
of the $\si$--bialgebras $T(D)$ and $T(C)$ that extends $\psi_1\ty$.

Applying this result to $C = C(\eta)$ and $D = C(\eta)^{\rm op}$, we obtain
a unique dual pairing $\psi$ of the $\si$--bialgebras
$T(C(\eta))^{\rm cop} = T(C(\eta)^{\rm op})$ and $T(C(\eta))$ such that
 \[ \psi(X_{ik} \ty, X_{j\ell})
     = \si(\eta_i - \eta_k \ty, \eta_j - \eta_{\ell}) R_{ij,\ty k\ty\ell}
                               \quad\mbox{for all $i,j,k,\ell \in I$} \,. \]
Let us next investigate under which conditions this form induces a bilinear
form on $A(R)$. For this to be the case it is necessary that
 \[ \psi(X_{ij,\ty k\ty\ell}\tty, X_{ab}) = 0
               \quad\mbox{and}\quad
    \psi(X_{ab}\tty, X_{ij,\ty k\ty\ell}) = 0 \]
for all $i,j,k,\ell,a,b \in I$. Each of these conditions is satisfied if and
only if the graded Yang--Baxter equation holds for $R\ty\ty$. Conversely, if
this is the case, the bilinear form $\psi$ indeed induces a bilinear form
$\vr$ on $A(R)$. Obviously, $\vr$ is a dual pairing of the $\si$--bialgebras
$A(R)^{\rm cop}$ and $A(R)$.

Our next task is to show that the generalized commutation relations
\reff{comrel} are satisfied. The basic first step of the proof consists in
checking that these relations hold for $a = t_{ik}$ and $b = t_{j\ell}\ty\ty$,
for all $i,j,k,\ell \in I$. But up to unimportant overall factors, these are
just the defining relations for $A(R)$.

Finally, we have to prove that $\vr$ is convolution invertible. If $\vr'$ is
an inverse of $\vr\ty\ty$, the Eqns.~\reff{rhoinv}, applied to $a = t_{ik}$ and
$b = t_{j\ell}\ty\ty$, show that $R$ is invertible: In fact, the inverse of
$R$ is given in terms of $\vr'$ exactly as $R$ is given in terms of
$\vr\tty$ (see Eqn.~\reff{rhott}). Accordingly, we now assume that $R$ is
invertible.

It is well--known (and easy to see) that if $R$ satisfies the graded
Yang--Baxter equation, then so does $\tilde{R} = P R^{-1\!} P$. Moreover,
$\tilde{R}$ and $R$ define the same $\si$--bialgebra: We have
 \be A(\tilde{R}) = A(R) \;. \label{ARti} \ee
Consequently, we can apply the foregoing results with $R$ replaced by
$\tilde{R}\ty\ty$, which gives another dual pairing $\tilde{\rho}$ of the
$\si$--bialgebras $A(R)^{\rm cop}$ and $A(R)$. According to the preceding
discussion, we expect that $\vr' = \tilde{\rho} \comp P$ is the inverse
of $\vr\ty\ty$, and this is indeed the case.

\sect{Graded comodules and comodule algebras \\ over $A(R)$ \vspace{-1ex}}
We keep the notation introduced at the beginning of Section 2 and in
Section 3. In the present section we are going to discuss certain natural
graded $A(R)$--comodules and $A(R)$--comodule algebras. Actually, we make a
slight generalization. The point is the following. At last, we are not so
much interested in the $\si$--bialgebras $A(R)$ themselves, but rather in
certain $\si$--Hopf algebra envelopes therefrom. Typically, the latter are
obtained from $A(R)$ by adjoining the inverses of certain elements of $A(R)$,
or by requiring additional relations (i.e., by going to a quotient). In any 
case, such an envelope will be a $\si$--bialgebra $A\tty$, and there will be
a natural homomorphism of $\si$--bialgebras $A(R) \rar A\tty$. For
convenience, the natural image of $t_{ij}$ in $A$ will also be denoted by
$t_{ij}\tty$. This will be the setting considered in the present section.
Stated differently, $A$ will be a $\si$--bialgebra containing certain
elements $t_{ij}\tty$; $i,j \in I$, such that $t_{ij}$ is homogeneous of
degree $\eta_j - \eta_i\tty$, for all $i,j \in I$, such that the relation
\reff{RTTA} (or its equivalents \reff{RTTAex}, \reff{RTTbr}) is satisfied,
and such that the Eqns.~\reff{delT} and \reff{epsT} hold true.

To begin with, we note that $V$ is a graded right $A$--comodule in a natural
way: The structure map
 \[ \de_1 : V \lra V \ot A \]
is defined by
 \be \de_1(e_i) = \sum_{j \in I} \si(\eta_i \tty, \eta_j - \eta_i) \tty
                                       e_j \ot t_{ji} \;, \label{rAcom} \ee
for all $i \in I$ (see Eqn.~\reff{rcom}).

Similarly, $\Vag$ is a graded left $A$--comodule: The structure map
 \[ \de'_1 : \Vag \lra A \ot \Vag \]
is defined by
 \be \de'_1(e'_i) = \sum_{j \in I} \si(\eta_j \tty, \eta_i - \eta_j) \tty
                                      t_{ij} \ot e'_j \;, \label{lAcom} \ee
for all $i \in I$. We note that the graded comodule $\Vag$ is dual to the
graded comodule $V$ in the sense that 
 \be (\id_A \ot \la \;\; \ra) \comp (\de'_1 \ot \id_V)
     \,=\, (\la \;\; \ra \ot \id_A) \comp (\id_{\Vag} \ot \de_1) \;.
                                             \label{dual} \vspace{0.8ex} \ee
(These are maps of $\Vag \ot V$ into $A\tty$, and
$\la \;\; \ra : \Vag \ot V \rar \KK$ is given by the dual pairing.)
Obviously, if one of the structure maps $\de_1$ and $\de_1'$ is known, the
other is uniquely fixed by this equation.
\vspace{2.0ex}

\noindent
{\em Remark 4.1.} At this point we insert a remark that could have been made
earlier. Looking at the Eqns.~\reff{rAcom} and \reff{lAcom}, we are tempted
to introduce the elements
 \be \tti_{ij} = \si(\eta_j \ty, \eta_i - \eta_j) \ty\ty t_{ij}
                       \quad\mbox{for all $i,j \in I$} \,. \label{tmod} \ee
Then these equations take the familiar form
 \[ \de_1(e_i) = \sum_{j \in I} e_j \ot \tti_{ji} \quad,\quad
                          \de'_1(e'_i) = \sum_{j \in I} \tti_{ij} \ot e'_j \]
for all $i \in I$, which implies that
 \[ \De(\tti_{ij}) = \sum_{k \in I} \tti_{ik} \ot \tti_{kj} \]
for all $i,j \in I$. (Needless to say, this equation can easily be checked
directly.) Summarizing, the elements $\tti_{ij}$ seem to have properties
which are closer to what we are accustomed to from the non--graded case.
However, unfortunately I do not know how to write the commutation relations
for the $\tti_{ij}$ in a factorized form like Eqn.~\reff{RTTA}. My conclusion
is that, depending on the problem at hand, both types of generators may have
advantages. (Actually, in the proof of Theorem 1 I have mostly used the
generators $\tti_{ij}\tty$, whereas for the proof of Eqn.~\reff{ARti} the
factorization in Eqn.~\reff{RTTA} is of great help.) \vspace{2.0ex}

Once the graded $A$--comodule structures of $V$ and $\Vag$ are defined, they
can be used to convert $T(V)$ and $T(\Vag)$ into graded $A$--comodule
algebras. Thus, there exists a unique homomorphism of graded algebras
 \[ \de : T(V) \lra T(V) \otb A \]
such that Eqn.~\reff{rAcom} (with $\de_1$ replaced by $\de\tty$) is satisfied,
and it is easy to see that this map converts $T(V)$ into a graded right
$A$--comodule, hence into a graded right $A$--comodule algebra. Note that,
for every integer $n \geq 0\ty\ty$, the subspace $T_n(V) = V^{\ot n}$ is a
graded subcomodule of $T(V)$, in fact, this is just the $n\tty$th tensorial
power of the graded comodule $V$. Its structure map is the map
 \be \de_n : T_n(V) \lra T_n(V) \ot A \label{deln} \ee
induced by $\de\tty$. (The map $\de_1$ is the map introduced in
Eqn.~\reff{rAcom}.)

Similarly, there exists a unique homomorphism of graded algebras
 \[ \de' : T(\Vag) \lra A \otb T(\Vag) \]
such that Eqn.~\reff{lAcom} (with $\de'_1$ replaced by $\de'\tty$) is
satisfied, and it is easy to see that this map converts $T(\Vag)$ into a
graded left $A$--comodule algebra. For every integer $n \geq 0\ty\ty$, the
subspace $T_n(\Vag) = (\Vag)^{\ot n}$ is a graded subcomodule of $T(\Vag)$,
in fact, this is just the $n\tty$th tensorial power of the graded comodule
$\Vag$. Its structure map is the map
 \be \de'_n : T_n(\Vag) \lra A \ot T_n(\Vag) \label{delprn} \ee
induced by $\de'$. (The map $\de'_1$ is the map introduced in
Eqn.~\reff{lAcom}.)

The graded comodule $T_n(\Vag)$ is dual to the graded comodule $T_n(V)$.
Indeed, since the case $n = 0$ is trivial, we may assume that $n \geq 1\ty$.
Let us define a canonical pairing
 \[ \la \;,\: \ra : T_n(\Vag) \ti T_n(V) \lra \KK \vspace{-1.0ex} \]
by
 \[ \la x'_1 \ot \ldots \ot x'_n \tty, x_1 \ot \ldots \ot x_n \ra
\,=\, \prod_{1 \leq a < b \leq n} \!\!\!\! \si(\xi'_b \tty,\ty \xi_a)
      \prod_{1 \leq c \leq n} \! \la x'_c \tty,\ty x_c \ra \]
for all $x_1 \ty,\ldots, x_n \in V$; $x'_1 \ty,\ldots, x'_n \in \Vag$,
which are homogeneous of the degrees $\xi_1 \ty,\ldots, \xi_n \ty,
\xi'_1 \ty,\ldots, \xi'_n\tty$, respectively. Then $\de_n$ and $\de'_n$
are related by an equation analogous to Eqn.~\reff{dual}.

The algebra $T(V)$ is free, with free generators $e_i$\,; $i \in I$. In
order to obtain more interesting examples of graded $A$--comodule algebras,
we have to construct suitable quotients of $T(V)$. In principle, this is
easily done (see Ref.~\cite{MMo}). Let $W$ be any graded $A$--subcomodule of
$T(V)$, let $J(W)$ be the ideal of $T(V)$ that is generated by $W$, and let
 \[ B_r(W) = T(V)/J(W) \]
be the corresponding quotient. Of course, $B_r(W)$ is a graded algebra. By
assumption, we have
 \[ \de(W) \subset W \ot A \;. \]
Obviously, this implies that
 \[ \de(J(W)) \subset J(W) \ot A \;. \]
Consequently, the map $\de$ induces a map
 \[ \de : B_r(W) \lra B_r(W) \otb A \]
(note the abuse of notation), and it is clear that this map converts $B_r(W)$
into a graded right $A$--comodule algebra. Let us denote the canonical image
of $e_i$ in $B_r(W)$ by $x_i\ty\ty$, for all $i \in I$. Then $B_r(W)$ is the
universal graded algebra, generated by elements $x_i$\,; $i \in I$, which
are homogeneous of degree $\eta_i$ and satisfy the relations ``given by
$W$\,''. Moreover, the comodule structure of $B_r(W)$ is fixed by
 \be \de(x_i) = \sum_{j \in I} \si(\eta_i \tty, \eta_j - \eta_i) \tty
           x_j \ot t_{ji} \quad\mbox{for all $i \in I$} \,. \label{delx} \ee

Instead of writing the defining relations in general, let us give a simple
example. Suppose that $W$ is a graded subcomodule of $T_2(V) = V \ot V$. Let
$(a^s)_{1 \leq s \leq r}$ be a family of elements of $W$ that generates the
vector space $W$, and let us write
 \[ a^s = \sum_{i,\ty j \in I}\! a^s_{ij} \ty\ty e_i \ot e_j
                                    \quad\mbox{for $1 \leq s \leq r$} \;. \]
Then the defining relations of $B_r(W)$ can be written in the form
 \[ \sum_{i,\ty j \in I}\! a^s_{ij} \ty\ty x_i \ty\ty x_j = 0
                                    \quad\mbox{for $1 \leq s \leq r$} \;. \]
 
In practice, it might be difficult to prove that a graded subspace of $T(V)$
is a graded subcomodule. Indeed, all we know about the algebra structure of
$A$ are the commutation relations \reff{RTTAex}. It is at this point where
the techniques of Ref.~\cite{RTF} are of great help.

Let $S \in \Lgr(V \ot V)$ be homogeneous of degree zero. As usual, we write
 \be S = \sum_{ijk\ty\ell \in I} \!\!
         S_{ij,\ty k\ty\ell}\ty\ty E_{ik} \otb E_{j\ell} \;. \label{Snot} \ee
We recall that the homogeneity condition is equivalent to the requirement that
 \[ S_{ij,\ty k\ty\ell} = 0
           \quad\mbox{if $\eta_i + \eta_j \neq \eta_k + \eta_{\ell}$} \;. \]
Let us also introduce the elements
 \[ T_1\ty,T_2 \in \Lgr(V) \otb \Lgr(V) \otb A \simeq \Lgr(V \ot V) \otb A \]
as in Eqn.~\reff{ToTt}. Then $S$ is an endomorphism of the graded
$A$--comodule $V \ot V$ if and only if
  \be (S \ot 1) T_1 T_2 = T_1 \ty T_2 (S \ot 1) \;. \label{STT} \ee
Explicitly, this relation is equivalent to
 \be \sum_{a,b \in I} \si(\eta_j - \eta_{\ell} \ty, \eta_a - \eta_k)
                   \tty S_{ij,\ty ab} \ty\ty t_{ak} \ty\ty t_{b\ty\ell}
   \,=\, \sum_{a,b \in I} \si(\eta_j - \eta_b \ty, \eta_i - \eta_k)
           \ty\ty t_{ia} \ty\ty t_{jb}\ty\ty S_{ab,\ty k\ty\ell} \;,
                                                          \label{STTex} \ee
for all $i,j,k,\ell \in I$.

A look at Eqn.~\reff{RTTbr} shows that $S = \Rh$ satisfies Eqn.~\reff{STT}.
Indeed, the algebra $A(R)$ is {\em defined} such that $\Rh$ is an endomorphism
\vspace{0.2ex}
of the graded $A(R)$--comodule $V \ot V$. It follows that every polynomial in
$\Rh$ satisfies this equation as well. Since the image and the kernel of
every such $S$ are graded $A$--subcomodules of $V \ot V$ (and hence of $T(V)$),
we have found some natural candidates for the subcomodule $W$. 

We can now comment on how the present formulation generalizes the approach of
Ref.~\cite{RTF}. In that reference, the $A(R)$--comodule $V \ot V$ is
completely reducible (provided that $q$ is generic), and the subcomodules can
be constructed as the images of suitable polynomials in $\Rh\ty\ty$. As we
are going to see, for the quantum supergroups $\mbox{SPO}_q(2n \ty|\ty 2n)$
this is not true anymore, and we are forced to use the generalized method
described above.

The graded left $A$--comodule algebras can be constucted in a completely
analogous way, and hence it should not be necessary to give all the details.
In particular, for every graded subcomodule $W'$ of $T(\Vag)$, we obtain a
graded left $A$-comodule algebra $B_{\ell}(W')$. However, one point should
be mentioned. In this second case, we need to know the endomorphisms of the
graded $A$--comodule $\Vag \ot \Vag$. Every element of $\Lgr(\Vag \ot \Vag)$
can be written as the graded transpose $S^{\rm T}$ of some element
$S \in \Lgr(V \ot V)$, and $S^{\rm T}$ is an endomorphism of the graded
$A$--comodule $\Vag \ot \Vag$ if and only if $S$ is an endomorphism of the
graded $A$--comodule $V \ot V$, i.e., if and only if the equivalent conditions
\reff{STT} and \reff{STTex} are satisfied. Using the notation \reff{Snot} and
assuming that $S$ is homogeneous of degree zero, we have
 \[ S^{\rm T}(e'_i \ot e'_j)
       \,=\, \sum_{k,\tty\ell \in I} \! \si(\eta_j \ty, \eta_k - \eta_i)
                \tty\tty S_{ij,\ty k\ty\ell} \,\ty e'_k \ot e'_{\ell} \;. \]

\sect{The quantum supergroup $\SP$ \vspace{-1ex}}
In the present and in the subsequent section we are going to apply the
general theory to the $R$--matrix that has been calculated in Ref.~\cite{Srm},
namely, to the $R$--matrix of $\Uq$ in the vector representation.
Correspondingly, we shall use the notation introduced in that reference.
First of all, this implies that in the following we choose
 \[ \Ga = \ZZ_2 = \{\ty\ol{0},\ol{1}\ty\} \]
and define
 \[ \si(\al\ty,\bet) = (-1)^{\al\bet}
                             \quad\mbox{for all $\al,\bet \in \ZZ_2$} \;. \]
Moreover, we set
 \[ r = n + m \quad,\quad d = n - m \;, \]
and the index set $I$ is equal to
 \[ I \,=\, \{-r,-r+1,\ldots,-2\ty,-1,1,2\ty,\ldots,r-1,r\} \;. \]
The degrees $\eta_i\ty\ty$; $i \in I$, are given by
 \[ \eta_i \,=\, \left\{ \begin{array}{ll}
                           \ol{0} & \mbox{for $1 \leq |i| \leq n$} \\[0.5ex]
                           \ol{1} & \mbox{for $n+1 \leq |i| \leq r$} \;.
                         \end{array} \right. \vspace{-0.5ex} \]
We recall that
 \[ \si_i = \si(\eta_i \ty, \eta_i) \quad\mbox{and}\quad
                      \si_{i,\ty j} = \si(\eta_i \ty, \eta_j)
                                     \quad\mbox{for all $i,j \in I$} \,. \]
In the present setting, $V$ is the vector module of $\Uq$.
\vspace{0.2ex}
It carries a $\Uq$--invariant bilinear form $b^q$, whose matrix
\vspace{0.4ex}
$C^q = (C^q_{i,\ty j})$ with respect to the basis $(e_i)_{i \in I}$ has been
given at the end of Section 4 of Ref.~\cite{Srm}. Finally, we note that the
$R$--matrix $R$ and the braid generator $\Rh$ are given by
  \[ \begin{array}{rcl}
      R &\! = \!& \dis{\sum_i q^{\si_i} \ty E_{i,\ty i} \otb E_{i,\ty i}
              \,+\, \sum_i q^{-\si_i}
                                 \ty E_{i,\ty i} \otb E_{-i,-i}}   \\[2.5ex]
        &   &\hspace{-1.2em}
             {} + \dis{\sum_{i \neq j,-j} \!
                                     E_{i,\ty i} \otb E_{j,\ty j}} \\[2.5ex]
        &   &\hspace{-1.2em}
             {} + (q - q^{-1})\ty \dis{\sum_{i < j}
                        \si_i \ty\ty E_{j,\ty i} \otb E_{i,\ty j}} \\[2.5ex]
        &   &\hspace{-1.2em}
             {} - (q - q^{-1})\ty \dis{\sum_{i < j}
                                   \si_i \tty \si_j \tty \si_{i,\ty j} \tty
            \Cqi_{-j,\ty j} \tty\ty C^q_{i,-i} \, E_{j,\ty i} \otb E_{-j,-i}}
     \end{array} \]
  \[ \begin{array}{rcl}
    \Rh &\! = \!& \dis{\sum_i \si_i \ty\tty
                                  q^{\si_i} \ty E_{i,\ty i} \otb E_{i,\ty i}
              \,+\, \sum_i \si_i \tty\ty q^{-\si_i}
                                 \ty E_{-i,\ty i} \otb E_{i,-i}}   \\[2.5ex]
        &   &\hspace{-1.2em}
             {} + \dis{\sum_{i \neq j,-j} \! \si_i \tty
                                     E_{j,\ty i} \otb E_{i,\ty j}} \\[2.5ex]
        &   &\hspace{-1.2em}
             {} + (q - q^{-1})\ty \dis{\sum_{i < j}
                                     E_{i,\ty i} \otb E_{j,\ty j}} \\[2.5ex]
        &   &\hspace{-1.2em}
         {} - (q - q^{-1})\ty \dis{\sum_{i < j} \si_i \tty \si_{i,\ty j} \tty
         \Cqi_{-j,\ty j} \tty\ty C^q_{i,-i} \, E_{-j,\ty i} \otb E_{j,-i}} \;,
     \end{array} \]
and that the operator $K$ in $V \ot V$ is equal to
  \[ K \,=\, \sum_{i,\ty j,\ty k,\ty\ell} \si_{k,\ty j} \tty \si_{k,\ty\ell}
             \ty\ty \Cqi_{i,\ty j} \ty\ty C^q_{k,\ty\ell}
                                   \, E_{i,\ty k} \otb E_{j,\ty\ell} \;. \]
\noindent
{\em Remark 5.1.} In Ref.~\cite{Srm} we have assumed that the base field
$\KK$ is equal to $\CC\,$ (since this assumption has been made by Yamane in
Ref.~\cite{Yam}). Obviously, the $R$--matrix obtained in Ref.~\cite{Srm}
makes sense under the present more general assumptions, and the properties
derived there are still true. Actually, also the definition of the quantum
superalgebra $\Uq$, its properties, and the discussion of the vector module
and of its tensorial square all remain valid. (Indeed, according to my
understanding Yamane's entire paper is correct with $\CC$ replaced by
$\KK\,$.) For convenience, we still assume that $q$ is not a root of unity,
although this condition could be relaxed. \vspace{2.0ex}

The $\si$--bialgebras $A(R)$ don't have an antipode (except for the
degenerate case where $I$ is empty): This follows immediately by means of the
$\ZZ$--gradation of these algebras. Hence some amendments are necessary if we
want to obtain a $\si$--Hopf algebra. In view of the comments made in
Section 2, it is easy to guess how this can be achieved under the present
assumptions. 

Let $T^{\rm\ty st}$ be the super--transpose of the matrix $T$, defined by
 \[ (T^{\rm\ty st})_{k,\ty\ell}
           = \si(\eta_{\ell}\ty, \eta_{\ell} - \eta_k) \ty\ty t_{\ell,k}
                                   \quad\mbox{for all $k,\ell \in I$} \,, \]
(see Eqn.~\reff{piT}), and let $F_r$ be the $I \ti I$--matrix over $\KK$\,,
defined by
 \[ (F_r)_{i,\ty j} = \si_{j,\ty i} \ty\ty C^q_{i,\ty j}
                       = \si_i \ty\ty C^q_{i,\ty j}
                          = C^q_{i,\ty j} \tty\ty \si_j 
                                          \quad\mbox{for all $i,j \in I$} \]
(see Eqn.~\reff{Fr}). Then the Eqns.~\reff{Spipi}, \reff{piS}, and \reff{Scpi}
suggest to define the $I \ti I$--matrix (with elements in $A(R)$)
 \[ T' = F_r^{-1}\tty T^{\rm\ty st} F_r \;, \]
explicitly:
  \[ t'_{i,\ty j} \,=\, \si_i \ty\ty \si_{i,\ty j} \ty\ty
                                ((C^q)^{-1})_{i,-i} \tty\ty
                                     t_{-j,-i} \tty\ty C^q_{-j,\ty j} \;,
                                                           \vspace{0.8ex} \]
and to investigate, for all $i,j \in I$, the expressions
  \[ \begin{array}{rcl}
     Q'_{i,\ty j}
        \eqq \dis{\sum_{k \in I}
                       \si(\eta_k - \eta_i \ty, \eta_k - \eta_j)\ty\ty
                                  t'_{i,\ty k} \ty\ty t_{k,\ty j} } \\[2.5ex]
    Q''_{i,\ty j}
        \eqq \dis{\sum_{k \in I}
                       \si(\eta_k - \eta_i \ty, \eta_k - \eta_j)\ty\ty
                                      t_{i,\ty k} \ty\ty t'_{k,\ty j} } \;.
     \end{array} \vspace{-0.8ex} \]

This can be carried out as in Ref.~\cite{KSc}. The basic observation is
\vspace{0.2ex}
that the operator $K$ can be written as a polynomial in $\Rh\ty\ty$ (see
Eqn.~(7.3) of Ref.~\cite{Srm}). This implies that the relations \reff{STTex}
are satisfied for $S = K\ty\ty$, and it is easy to evaluate these relations.
                                                                    \\[1.5ex]
If $i \neq -j$ and $k \neq -\ell\,$, both sides of these relations are
trivially equal to zero, and no information can be obtained. \\[1.5ex]
If $i = -j$ and $k \neq -\ell\,$, it follows that
  \[ \sum_{a,\ty b \in I} \si_{a,\ty\ell} \tty\ty \si_{b,\ty a} \tty\ty
          C^q_{a,\ty b} \tty\ty t_{a,\ty k} \tty\ty t_{b,\ty\ell} \,=\, 0
                                    \quad\mbox{for $k \neq -\ell$} \;. \]
If $i \neq -j$ and $k = -\ell\,$, we find that
  \[ \sum_{a,\ty b \in I} \si_{i,\ty b}\tty\ty ((C^q)^{-1})_{a,\ty b} \tty\ty
                                  t_{i,\ty a} \tty\ty t_{j,\ty b} \,=\, 0
                                      \quad\mbox{for $i \neq -j$} \;.  \]
If $i = -j$ and $k = -\ell\tty$, we set
  \[ \begin{array}{rcl}
     Q' \eqq \dis{(C^q_{k,-k})^{-1} \sum_{a,\ty b \in I}
                   \si_{a,\ty k} \tty\ty \si_{b,\ty a} \tty\ty
                         C^q_{a,\ty b} \tty\ty
                                   t_{a,\ty k} \tty\ty t_{b,-k} } \\[3.0ex]
    Q'' \eqq \dis{(((C^q)^{-1})_{i,-i})^{-1} \ty \si_i
                   \sum_{a,\ty b \in I}
                       \si_{i,\ty b} \tty\ty ((C^q)^{-1})_{a,\ty b} \tty\ty
                                    t_{i,\ty a} \tty\ty t_{-i,\ty b} } \;.
     \end{array} \]
Then the relation under investigation says that
  \[ Q' = Q'' \;. \]
Obviously, $Q'$ does not depend on $i,j\ty$, and $Q''$ does not depend on
$k,\ell\tty$. This implies that $Q'$ and $Q''$ do not depend on
$i,j,k,\ell\tty$, and the equations above can be summarized as follows:
  \bea \dis{ \sum_{a,\ty b \in I}\!
                 \si_{\ell,\ty b} \tty\ty \si_{b,\ty a} \tty\ty
                         C^q_{a,\ty b} \tty\ty
                               t_{a,\ty k} \tty\ty t_{b,\ty\ell} }
         \eqq C^q_{k,\ty\ell}\, Q'                          \label{invb} \\
       \dis{ \sum_{a,\ty b \in I}\!
                 \si_{i,\ty b} \tty\ty ((C^q)^{-1})_{a,\ty b} \tty\ty
                              t_{i,\ty a} \tty\ty t_{j,\ty b} }
         \eqq \si_{i,\ty j} \ty\ty ((C^q)^{-1})_{i,\ty j}\, Q'' \;.
                                                            \label{inva} \eea
\noindent
{\em Remark 5.2.} Using the generators $\tti_{i,\ty j}$ as defined in
Eqn.~\reff{tmod}, i.e.,
 \[ \tti_{i,\ty j} = \si_j \ty\ty \si_{j,\ty i} \ty\ty t_{i,\ty j}
                                     \quad\mbox{for all $i,j \in I$} \,, \]
and defining the super--transpose $\tilde{T}^{\rm\ty st}$ of the
$I \ti I$--matrix $\tilde{T} = (\tti_{k,\ty\ell})$ by
 \[ (\tilde{T}^{\rm\ty st})_{k,\ty\ell}
          = \si_{\ell} \ty\ty \si_{k,\ty\ell} \ty\ty \tti_{\ell,\ty k}
                                      \quad\mbox{for all $k,\ell \in I$} \]
(see Eqn.~\reff{piT}), it is easy to see that the relations \reff{invb} and
\reff{inva} can be written in the form
 \[ \tilde{T}^{\rm\ty st} \ty C^q \ty\ty \tilde{T} = \ty Q' C^q \quad,\quad
    \tilde{T} \ty (C^q)^{\ty\!-1} \ty\ty \tilde{T}^{\rm\ty st}
                            = \ty Q'' (C^q)^{\ty\!-1} \,. \vspace{1.0ex} \]

Now it is easy to see that
  \be  Q'_{i,\ty j} \,=\, \de_{i,\ty j} \ty\ty Q' \quad,\quad 
     Q''_{i,\ty j} \,=\, \de_{i,\ty j} \ty\ty Q'' \;. \label{TTpr} \ee
In the following, we set
  \[ Q' = Q'' = Q \;. \] 
Obviously, the element $Q$ is homogeneous of degree zero. Using the
Eqns.~\reff{invb} and \reff{inva} once again, we can also show that $Q$ is
group--like, i.e., that
  \[ \De(Q) = Q \ot Q \quad,\quad \ve(Q) = 1 \;. \]

Finally, we prove that $Q$ lies in the center of $A(R)$. Consider the
following elements of the graded tensor product $\Lgr(V) \otb A(R)$\,:
 \[ T_0 = \sum_{i,\ty j \in I} E_{i,\ty j} \ot t_{i,\ty j}
                     \quad,\quad
   T'_0 = \sum_{i,\ty j \in I} E_{i,\ty j} \ot t'_{i,\ty j} \;. \]
Then the Eqns.~\reff{TTpr} can be written in the following form:
  \be T'_0 \ty\ty T_0 = T_0 \ty\ty T'_0 \ty = \ty \II \ot Q \label{Tinv} \ee
(where we have set $\II = \id_V$). These equations imply that
  \[ T_0 \tty (\ty\II \ot Q) = T_0 \ty\ty T'_0 \tty\ty T_0
                                            = (\ty\II \ot Q) \ty T_0 \;. \]
Since $Q$ is homogeneous of degree zero, this is equivalent to
  \[ t_{i,\ty j} \, Q = Q \, t_{i,\ty j}
                                      \quad\mbox{for all $i,j \in I$} \,, \]
which proves our claim.

After these preparations, we can define the quantum supergroup $\SP$, as
follows. Let $J(Q)$ be the ideal of $A(R)$ that is generated by $Q - 1\tty$.
Since $Q$ is homogeneous of degree zero and group--like, $J(Q)$ is a graded
biideal of $A(R)$. Consequently, the quotient
 \be \SP = A(R)/J(Q) \label{defSP} \ee
inherits from $A(R)$ the structure of a bi--superalgebra (i.e., a
$\si$--bialgebra with respect to the commutation factor $\si$ of
supersymmetry). Let
 \[ \om : A(R) \lra \SP \]
be the canonical map. Until further notice, we are going to use the notation
 \[ \om(t_{i,\ty j}) = \tb_{i,\ty j} \quad\mbox{for all $i,j \in I$} \,. \]
Similarly, we introduce the elements $\tb\ty'_{i,\ty j}$ and the
$I \ti I$--matrices $\ol{T}$ and $\ol{T}\ty'$. 

Let us now show that $\SP$ is a Hopf superalgebra. According to the preceding
results we expect that the antipode $S$ of $\SP$ is given by
 \be S(\ol{T}\ty) = \ol{T}\ty' \label{anti} \ee
(to be interpreted element--wise). Consequently, we have to show that there
exists a homomorphism of graded algebras
 \be S : \SP \lra \SP^{\rm aop} \label{Shom} \ee
such that Eqn.~\reff{anti} is satisfied.

To prove this, we work in the graded tensor product
 \[ \Lgr(V) \otb \Lgr(V) \otb \SP \simeq \Lgr(V \ot V) \otb \SP \;. \]
In that algebra, we have (with our usual notation)
 \[ (R \ot 1) \ol{T}_{\!1} \ty \ol{T}_{\!2}
                        = \ol{T}_{\!2} \ty\ty \ol{T}_{\!1} (R \ot 1) \;, \]
and also
 \[ \ol{T}_{\!1}^{\ty\ty-1} = \ol{T}\ty'_{\!1} \quad,\quad
                               \ol{T}_{\!2}^{\ty\ty-1} = \ol{T}\ty'_{\!2} \] 
(see Eqn.~\reff{Tinv}). This implies that
 \[ (R \ot 1) \ol{T}\ty'_{\!2} \tty \ol{T}\ty'_{\!1}
               = \ol{T}\ty'_{\!1} \ty\tty \ol{T}\ty'_{\!2} (R \ot 1) \;. \]
A moment's thought then shows that this is the RTT--{\tty}relation \reff{RTTA}
for $\ol{T}\ty'$ in $\SP^{\rm aop}$.

This result implies that there exists a homomorphism of graded algebras
 \[ S_0 : A(R) \lra \SP^{\rm aop} \]
such that
 \[ S_0(T) = \ol{T}\ty' \;. \]
It is easy to check that
 \[ S_0(Q') = \om(Q'') = 1 \quad,\quad S_0(Q'') = \om(Q') = 1 \;. \]
Consequently, $S_0$ induces the homomorphism \reff{Shom} we are looking for.
It is now easy to see that $S$ is an antipode: The identities to be proved
hold (by construction) when evaluated on the generators $\tb_{i,\ty j}\tty$
of $\SP$, and this implies that they hold on all of $\SP$.

Summarizing, we have shown that $\SP$ is a Hopf superalgebra. In the sequel,
we shall {\em simplify the notation} and write $t_{i,\ty j}$ instead of
$\tb_{i,\ty j}\ty\ty$. Then the antipode is uniquely fixed by the equation
 \[ S(T) = T' \,, \]
where
 \[ T' = F_r^{-1}\tty T^{\rm\ty st} F_r \;, \]
or more explicitly
 \[ t'_{i,\ty j} \,=\, \si_i \ty\ty \si_{i,\ty j} \ty\ty
                                ((C^q)^{-1})_{i,-i} \tty\ty
                                     t_{-j,-i} \tty\ty C^q_{-j,\ty j}
                       \quad\mbox{for all $i,j \in I$} \,. \vspace{0.8ex} \]
We note that
 \[ S^2(t_{i,\ty j}) = d_i \ty\ty d_{-j} \ty\ty t_{i,\ty j}
                                      \quad\mbox{for all $i,j \in I$} \,, \]
where
 \[ d_i = \si_i \ty\ty \Cqi_{i,-i} \ty\ty C^q_{i,-i}
                                        \quad\mbox{for all $i \in I$} \,. \]
Obviously, we have
 \[ d_{-i} = d_i^{-1} \quad\mbox{for all $i \in I$} \,. \]

In $\SP$, the relations \reff{invb} and \reff{inva} take the form
  \bea \dis{ \sum_{a,\ty b \in I}\!
                 \si_{\ell,\ty b} \tty\ty \si_{b,\ty a} \tty\ty
                         C^q_{a,\ty b} \tty\ty
                               t_{a,\ty k} \tty\ty t_{b,\ty\ell} }
         \eqq C^q_{k,\ty\ell}                             \label{invbSP} \\
       \dis{ \sum_{a,\ty b \in I}\!
                 \si_{i,\ty b} \tty\ty ((C^q)^{-1})_{a,\ty b} \tty\ty
                              t_{i,\ty a} \tty\ty t_{j,\ty b} }
         \eqq \si_{i,\ty j} \ty\ty ((C^q)^{-1})_{i,\ty j}
                                                        \label{invaSP} \eea
for all $i,j,k,\ell \in I$. These relations have a simple interpretation.
Eqn.~\reff{invbSP} is equivalent to the fact that the linear form
$\bt^q : V \ot V \rar \KK$ associated to the bilinear form $b^q$ is a
homomorphism of graded right $\SP$--comodules, and also to the fact that the
element
 \be a' = \sum_{i,\ty j \in I} \!
                    \si_{j,\ty i} \ty\ty C^q_{i,\ty j} \, e'_i \ot e'_j \ee
is an invariant of the graded left $\SP$--comodule $\Vag \ot \Vag$, in the
well--known sense that
 \[ \de'_2(a') = 1 \ot a' \;. \]
Similarly, Eqn.~\reff{invaSP} is equivalent to the fact that the linear form
 \[ \bt'^q : \Vag \ot \Vag \lra \KK \vspace{-0.5ex} \]
defined by
 \[ \bt'^q(e'_i \ot e'_j) = \si_{i,\ty j} \ty\ty \Cqi_{i,\ty j}
                      \quad\mbox{for all $i,j \in I$} \,, \vspace{0.5ex} \]
is a homomorphism of graded left $\SP$--comodules, and also to the fact that
the element
 \be a = \sum_{i,\ty j \in I} ((C^q)^{-1})_{i,\ty j} \, e_i \ot e_j
                                                           \label{defa} \ee
is an invariant of the graded right $\SP$--comodule $V \ot V$, in the sense
that
 \[ \de_2(a) = a \ot 1 \;. \]
Of course, this is exactly what we should guess.

According to Section 2, we expect that the Hopf superalgebra $\SP$ is
coquasitriangular. More precisely, we expect that the universal $r$--form
$\vr$ on $A(R)$ (as specified in Theorem 1) induces a universal $r$--form on
$\SP$.

In order to show this we have to prove that $\vr$ vanishes on $J(Q) \ti A(R)$
and on $A(R) \ti J(Q)$ (where $J(Q)$ is the ideal of $A(R)$ used in the
definition of $\SP$, see Eqn.~\reff{defSP}). The basic ingredients of the
proof are the equations
 \[ \vr(C^q_{k,\ty\ell}\tty(Q' - 1),\ty t_{i,\ty j}) = 0 \quad,\quad
            \vr(t_{i,\ty j} \ty\ty, C^q_{k,\ty\ell}\tty(Q' - 1)) = 0 \;, \]    
which hold for all $i,j,k,\ell \in I$. These are equivalent to the equations
 \[ R^{\tty{\rm st}_1} \comp (f_r \ot \II) \comp R = f_r \ot \II
                \quad,\quad
    R^{\tty{\rm st}_2} \comp (\II \ot f_{\ell}) \comp R
                                                  = \II \ot f_{\ell} \;, \]
respectively, which have been proved in Ref.~\cite{Srm} (see Eqn.~(7.5) of
that reference). We recall that $f_{\ell}$ and $f_r$ are the linear maps of
$V$ into $\Vag$ given by
 \[ f_{\ell}(e_j) = \sum_{i \in I} \ty C^q_{j,\ty i} \, e'_i
                                     \quad,\quad
         f_r(e_j) = \sum_{i \in I} \si_{j,\ty i}
                                            \ty\ty C^q_{i,\ty j} \, e'_i \]
for all $j \in I$.

According to the preceding remarks, the bilinear form $\vr$ induces a bilinear
form $\ol\vr$ on $\SP$. Obviously, this form satisfies the
Eqns.~\reff{comrel}\,--\,\reff{multr}. In order to prove that $\ol\vr$ is
convolution invertible, one might think we have to show that the inverse
$\vr'$ of $\vr$ also induces a bilinear form on $\SP$. However, it is not
necessary to prove this independently: Since $\SP$ is a Hopf superalgebra,
it is known (and easy to see) that the bilinear form $\ol\vr\tty'$ on $\SP$,
defined by
 \[ \ol\vr\tty'(b \ty, b'\ty) = \ol\vr(S(b),\ty  b'\ty)
                                \quad\mbox{for all $b\ty,b' \in \SP$} \;, \]
is an inverse of $\ol\vr\ty\ty$.

Simplifying the notation, we have proved the following theorem. \vspace{1.0ex}

\noindent
{\bf Theorem 2.} The Hopf superalgebra $\SP$ is coquasitriangular. More
precisely, there exists a unique universal $r$--form $\vr$ on $\SP$ such that
 \be \vr(t_{ik} \tty, t_{j\ell})
     = \si(\eta_i - \eta_k \tty, \eta_j - \eta_{\ell}) R_{ij,\ty k\ty\ell}
              \quad\mbox{for all $i,j,k,\ell \in I$} \,. \label{rhottSP} \ee
(We are using the notation introduced in Eqn.~\reff{Rnot}. Since $R$ is
homogeneous of degree zero, the factor
$\si(\eta_i - \eta_k \ty, \eta_j - \eta_{\ell})$ can be replaced by
$\si_i \ty\ty \si_k\ty\ty$, and also by
\vspace{1.0ex}
$\si_j \ty\ty \si_{\ell}\tty$.)

Finally (and once again according to Section 2) we expect that the Hopf
superalgebras $\Uq$ and $\SP$ are dual to each other in some sense. More
precisely, we expect that there exists a Hopf superalgebra pairing of $\SP$
and $\Uq$. This can be proved, as follows.

Consider the vector module $V$ of $\Uq$ (see Ref.~\cite{Srm}), and let $\pi$
be the graded representation of $\Uq$ afforded by it. Define the linear forms
$\pi_{i,\ty j} \in \Uq^{\circ}$ by
 \[ \pi_{i,\ty j}(h) = \la e'_i \tty, \pi(h)\ty e_j \ra
                  \quad\mbox{for all $i,j \in I$ and all $h \in \Uq$} \,. \]
We know that $\Rh$ is an endomorphism of the graded $\Uq$--module $V\ot V$.
Defining the elements $\pi_1$ and $\pi_2$ of the superalgebra
 \be \Lgr(V) \otb \Lgr(V) \otb \Uq^{\circ}
                      \simeq \Lgr(V \ot V) \otb \Uq^{\circ} \label{LLUq} \ee
as in Eqn.~\reff{pionetwo}, it is easy to see that this is equivalent to the
fact that the equation
 \[ (\Rh \ot \ve) \tty \pi_1 \tty \pi_2
                                  = \pi_1 \ty\ty \pi_2 \tty (\Rh \ot \ve) \]
holds in the algebra \reff{LLUq}, where $\ve$ is the counit of $\Uq$, i.e.,
the unit element of $\Uq^{\circ}$.

By definition of $A(R)$, this implies that there exists a unique algebra
homomorphism
 \[ \chi_0 : A(R) \lra \Uq^{\circ} \]
such that \vspace{-0.5ex}
 \[ \chi_0(t_{i,\ty j}) = \pi_{i,\ty j}
                       \quad\mbox{for all $i,j \in I$} \,. \vspace{0.5ex} \]
Obviously, $\chi_0$ is a bi--superalgebra homomorphism, and we have
 \[ \chi_0(C^q_{k,\ty\ell}\, Q')
          \,=\, \sum_{a,\ty b \in I}
                 \si_{\ell,\ty b} \tty\ty \si_{b,\ty a} \tty\ty
                       C^q_{a,\ty b} \tty\ty
                             \pi_{a,\ty k} \tty\ty \pi_{b,\ty\ell} 
                                   \quad\mbox{for all $k,\ell \in I$} \,. \]
Since the bilinear form $b^q$ on $V$ is $\Uq$--invariant (equivalently, since
the associated linear form $\bt^q$ on $V \ot V$ is a homomorphism of graded 
$\Uq$--modules), the right hand side is equal to
$C^q_{k,\ty\ell}\, \ve \ty\ty$. Consequently, $\chi_0$ induces a
bi--superalgebra homomorphism
 \[ \chi : \SP \lra \Uq^{\circ} \]
such that \vspace{-0.5ex}
 \[ \chi(t_{i,\ty j}) = \pi_{i,\ty j}
                       \quad\mbox{for all $i,j \in I$} \,. \vspace{0.5ex} \]
Since $\SP$ and $\Uq^{\circ}$ are Hopf superalgebras, $\chi$ is known to be
a Hopf superalgebra homomorphism. This implies (indeed, is equivalent to the
fact) that the bilinear form
 \[ \vp : \SP \ti \Uq \lra \KK \]
defined by
 \be \vp(a \tty, h) = \la \chi(a),\tty h \ra
      \quad\mbox{for all $a \in \SP$ and all $h \in \Uq$} \label{defphi} \ee
is a Hopf superalgebra pairing.

Summarizing, we have proved the following theorem. \vspace{1.0ex}

\noindent
{\bf Theorem 3.} There exists a unique Hopf superalgebra pairing
 \[ \vp : \SP \ti \Uq \lra \KK \]
such that
 \[ \vp(t_{i,\ty j} \ty\ty,\ty\ty h) \,=\, \pi_{i,\ty j}\tty(h)
              \quad\mbox{for all $i,j \in I$ and all $h \in \Uq$} \,.
                                                           \vspace{1.0ex} \]

Let $A$ be the image of $\SP$ in $\Uq^{\circ}$ under the homomorphism
$\chi$. Obviously, $A$ is the subalgebra of $\Uq^{\circ}$ that is
generated by the elements $\pi_{i,\ty j}\ty\ty$; $i,j \in I$. Actually, $A$
is a sub--Hopf--superalgebra of $\Uq^{\circ}$ (see the remark below
Eqn.~\reff{piT}). By definition, $A$ is dense in $\Uq^{\ast{\rm gr}}$ if and
only if it separates the elements of $\Uq$.  Using the definition
\reff{defphi}, we see that the implication
 \be \vp(a \tty, h) = 0 \;\;\mbox{for all $a \in \SP$}
                        \quad \Longrightarrow \quad h = 0 \label{Adense} \ee
is true for all $h \in \Uq$ if and only if $A$ is dense in
$\Uq^{\ast{\rm gr}}$. On the other hand, the implication
 \be \vp(a \tty, h) = 0 \;\;\mbox{for all $h \in \Uq$}
                        \quad \Longrightarrow \quad a = 0 \label{chiinj} \ee
is true for all $a \in \SP$ if and only if $\chi$ is injective (and hence
induces a Hopf superalgebra isomorphism of $\SP$ onto $A$). Unfortunately,
at present I am not able to prove or disprove one or even both of the
implications above. Actually, there is a simple reason to expect that
\reff{chiinj} is not correct (see Section 7).

\sect{Graded comodule algebras over $\SP$ \vspace{-1ex}}
We keep the notation of the preceding sections, in particular, that of
Section 5. In the following, we are going to apply the general results of
Section 4 to the case of the Hopf superalgebra $\SP$.

To begin with, we recall some well--known facts (actually, these are true
in the general framework of $\si$--bialgebras). If $W$ is a graded right
$\SP$--comod\-ule, the dual pairing of Theorem 3 can be used to convert $W$
into a graded left $\Uq$--module $\Wm$. The procedure is the same as that in
Section 2. Let
 \[ \de : W \lra W \ot \SP \]
be the structure map of $W$, and let $h \in \Uq$ and $x \in W$ be homogeneous
elements, of degrees $\eta$ and $\xi\ty\ty$, respectively. If
 \[ \de(x) = \sum_\su x^0_s \ot x^1_s \;, \vso \]
with $x^0_s \in W$ and $x^1_s \in \SP$, the action of $h$ on $x$ is defined by
 \[ h \cd x  = \si(\eta,\xi) \sum_\su x^0_s \, \vp(x^1_s \ty, h) \;. \vso \]
If $U$ is a second graded right $\SP$--comodule, and if $f \!:\!\ty U \rar W$
is a homomorphism of graded $\SP$--comodules, then
$f \!:\! U^{\rm mod} \rar \Wm$ is a homomorphism of graded $\Uq$--modules. In
particular, any graded $\SP$--subcomodule of $W$ is a graded $\Uq$--submodule
of $\Wm$. Also, if $W$ is the graded tensor product of the graded right
$\SP$--comodules $W_1\ty$,\,\ldots,\ty $W_p\ty\ty$, then $\Wm$ is the graded
tensor product of the graded left $\Uq$--modules $\Wm_1$,\,\ldots,\ty $\Wm_p$.
In addition, as the reader will expect, if $W$ is the vector comodule $V$ of
$\SP$, then $\Wm$ is the vector module $V$ of $\Uq$. Finally, if $B$ is a
graded right $\SP$--comodule algebra, the construction above converts $B$
into a graded left $\Uq$--module algebra (i.e., the structure maps of the
algebra $B$ are homomorphisms of graded $\Uq$--modules).

There is an analogous construction that converts graded left $\SP$--comodules
into graded right $\Uq$--modules, and its properties are completely analogous
to those mentioned above. We leave it to the reader to spell out the details.

After these preliminaries, we are ready to discuss the graded right
$\SP$--comodule algebras. Guided by Section 4, we should look for graded
$\SP$--subcomodules of $T(V)$. We restrict our attention to quadratic comodule
algebras. Correspondingly, we shall first determine the graded
$\SP$--subcomodules of
 \[ T_{\leq 2}(V) = T_0(V) \op  T_1(V) \op  T_2(V) \;. \] 
In view of the introductory remarks, we can do this in two steps: First, we
determine the graded $\Uq$--submodules of $T_{\leq 2}(V)$, then we try to
show that the $\Uq$--submodules we have found are, in fact, graded
$\SP$--subcomodules.

In the following, the $\SP$--comodules and $\Uq$--modules will simply be called
comodules and modules, respectively.

It is easy to see that every graded submodule of $T_{\leq 2}(V)$ is equal to
$U$ or $U \op T_1(V)$, where $U$ is a graded submodule of $T_0(V) \op T_2(V)$.
(Of course, we are not interested in the second case, since it leads to trivial
comodule algebras.) Similarly, every graded submodule of $T_0(V) \op T_2(V)$
is equal to $W$ or $W \op (V \ot V)_s\ty\ty$, where $W$ is a  graded submodule
of $T_0(V) \op (V \ot V)_a$ (we are using the notation of Ref.~\cite{Srm}).

Finally, every graded submodule $W$ of $T_0(V) \op (V \ot V)_a$ belongs to
one of the following two classes. \\[1.5ex]
1) The module $W$ is equal to $W_a$ or $T_0(V) \op W_a\ty\ty$, where $W_a$ is
a graded submodule of $(V \ot V)_a\ty\ty$. \\[1.5ex]
2) We have
 \[ W = \KK\ty (g - c) \op (V \ot V)^0_a \;, \]
where $g$ is any element of $(V \ot V)_a$ that does not belong to
$(V \ot V)^0_a\ty\ty$, and where $c \in \KK = T_0(V)$ is different from zero.

Our next task consists in checking which of these submodules are, in fact,
graded subcomodules of $T(V)$. This is done by means of the following
observations. (Unfortunately, in the case $n = m = 1\ty$, I haven't solved
this problem completely, see below.)

\noindent
Obviously, $T_1(V)$ is a graded subcomodule of $T(V)$.

Secondly, we know (see Ref.~\cite{Srm}) that the projector of $V \ot V$ onto
$(V \ot V)_s$ with kernel $(V \ot V)_a$ can be written as a polynomial in
$\Rh\ty\ty$. Since $\Rh$ (and hence the projector) is a homomorphism of
graded comodules, this implies that $(V \ot V)_s$ and $(V \ot V)_a$ are
graded subcomodules.

We know that $\bt^q : V \ot V \lra \KK$ is a homomorphism of graded comodules
(see Section 5), and that $(V \ot V)^0_a$ is the kernel of the restriction of
\vspace{0.2ex}
$\bt^q$ onto $(V \ot V)_a$ (see Ref.~\cite{Srm}). Hence $(V \ot V)^0_a$ is a
graded subcomodule as well. \vspace{2.0ex}

\noindent
{\em Remark 6.1.} We note that in the case $n = m$ there does not exist a
projector of $V \ot V$ onto $(V \ot V)^0_a$ that would be a homomorphism of
graded comodules. Hence in this case the projector method of Ref.~\cite{RTF}
cannot be applied. \vspace{2.0ex}

In Section 5 we have also noted that the element $a$ (see Eqn.~\reff{defa})
is an invariant of the comodule $V \ot V$. This shows that $\KK\,a$ is a
subcomodule of $V \ot V$.

In the case $n = m = 1\ty$, we also have to check whether $V_4$ and $\Vf$ are
graded subcomodules. Unfortunately, I haven't been able to solve this problem,
but I think there are good reasons to conjecture that they are not (see
Section 7).

Finally, each of the graded comodules of case 2) above is a graded subcomodule
as well. In fact, it is equal to the kernel of the linear form
 \[ f : T_0(V) \op (V \ot V)_a \lra \KK \]
defined by
 \[ f(\lam + u) = \bt^q(g)\ty\lam + c\,\bt^q(u)
           \quad\mbox{for all $\lam \in T_0(V)$, $u \in (V \ot V)_a$} \;, \]
which is a homomorphism of graded comodules since $\bt^q$ is such (recall
that $\bt^q(g)$ is not equal to zero).

Since, for each of the graded subcomodules of $T_{\leq 2}(V)$ found above,
a basis has been given in Ref.~\cite{Srm}, the defining relations of the
corresponding graded right $\SP$--comodule algebra $B_r$ can immediately be
written down.

In the following, we shall concentrate on the case 2) above, but allow for
$c = 0$ (i.e., we include the case $W = (V \ot V)_a\ty$). Note that, for
$c \neq 0\ty\ty$, the subcomodule $W$ only depends on the residue class
$c^{-1} g + (V \ot V)^0_a\ty\ty$. We fix our conventions by choosing
 \[ W = \KK\ty (t - q^{2d} c) \op (V \ot V)^0_a \;. \]
Here, $t \in (V \ot V)_a$ is the element specified in Eqn.~(5.17) of
Ref.~\cite{Srm}, and $c \in \KK$ is an arbitrary scalar. For convenience, we
have included the factor $q^{2d}\tty$; recall that $d = n - m\ty\ty$. The
corresponding graded right $\SP$--comodule algebra $B_r(W)$ will be denoted
by $\Arq$. For $c = 0$ it might be called the superalgebra of functions on
the symplecto--orthogonal quantum superspace, for $c \neq 0$ it is the quantum
Weyl superalgebra we wanted to construct (see Remark 6.2).

Using the basis of $(V \ot V)^0_a$ given in Ref.~\cite{Srm} and the formula
for $t$ given there, this comodule superalgebra can be described in terms of
generators and relations, as follows.

By definition, $\Arq$ is the universal superalgebra, generated by elements
$x_i\ty\ty$; $i \in I$, which are homogeneous of degree $\eta_i$ and satisfy
the following relations: \\[1.0ex]
The basis elements $a_{i,\ty i}$ with $\eta_i = \bar{1}$ give
  \[ x_i^2 = 0 \quad\mbox{for $\eta_i = \bar{1}$} \;, \]
the basis elements $a_{i,\ty j}$ with $i < j$ but $i \neq -j$ give
  \[ x_i \tty x_j = \si_{i,\ty j} \tty\ty q \ty\ty x_j \tty x_i 
                         \quad\mbox{for $i < j \ty,\: i \neq -j$} \;, \]
the tensors $a_j\ty\ty$, $2 \leq j \leq r\tty$, give
  \[ q^{-\si_{j-1}} x_{-j+1} \ty\ty x_{j-1}
                      - \si_{j-1} \ty\ty q \tty\ty x_{j-1} \ty\ty x_{-j+1}
     \,=\, \si_{j-1} \ty\ty \si_j \ty\ty x_{-j} \ty\ty x_j
         - \si_{j-1} \ty\ty q \tty\ty q^{-\si_j} \ty x_j \ty\ty x_{-j} \;, \]
and the tensor $t - q^{2d} c$ gives
  \[ x_{-1} \ty\ty x_1 - q^2 \tty x_1 \ty\ty x_{-1}
     \,=\, q^{2d} c - (q - q^{-1}) \sum_{i=2}^{n+m} \ty
                                 \Cqi_{i,-i} \tty\ty x_i \ty\ty x_{-i} \;.\]
The last $(r - 1) + 1$ relations are equivalent to the system
  \[ \begin{array}{rrcl}
    \!\!\mbox{(I)}\!
    & q^{\si_i} \ty x_{-i} \ty\ty x_i
			     - \si_i \ty\ty q^{-1} \ty x_i \ty\ty x_{-i}
    \eqq -(q - q^{-1}) \tty\tty \si_i \tty\tty C^q_{i,-i} \ty\ty
     \dis{\sum_{j < i} \tty ((C^q)^{-1})_{-j,\ty j} \tty\ty x_{-j} \ty\ty x_j
                    + \si_i \ty\ty C^q_{i,-i} \tty q^{2d} \ty c } \;,
     \end{array} \]
which, in turn, can be shown to be equivalent to the system
  \[ \begin{array}{rrcl}
    \!\!\!\mbox{(II)}\;\,
    & q^{-\si_i} x_i \ty\ty x_{-i} - \si_i \ty\ty q \ty\ty x_{-i} \ty\ty x_i
    \eqq (q - q^{-1}) \tty\tty \si_i \tty\tty C^q_{-i,i} \ty\ty
         \dis{\sum_{j < i} \tty ((C^q)^{-1})_{j,-j} \tty\ty x_j \ty\ty x_{-j}
                      + \si_i \ty\ty C^q_{-i,\ty i} \tty\ty c } \,. \;\,
     \end{array} \]
Both systems hold for $-r \leq i \leq -1\ty$ (and the summation is over all
$j \in I$ such that $j < i\ty$).

Finally, the structure map
 \[ \de : \Arq \lra \Arq \otb \SP \]
is fixed by Eqn.~\reff{delx}. 

It can be shown that the canonical image of the invariant $a \in V \ot V$ is
given by
 \be \sum_{j \in I} \tty ((C^q)^{-1})_{j,-j} \tty\ty x_j \ty\ty x_{-j}
         \,=\, (q^d - q^{-d}\tty)(q - q^{-1})^{-1} \ty\ty q^d \tty c \;.  
                                          \vspace{-1.5ex} \label{canima} \ee
For $n = m\ty\ty$, this element is equal to zero, as it should (since in this
case, we have $a \in (V \ot V)^0_a\tty$).

Using the braid generator $\Rh\ty\ty$, the defining relations of $\Arq$ can
be written in the following form. Let $M : A^r_q \ot A^r_q \rar A^r_q$ denote
the multiplication map of $\Arq$. Then the defining relations are equivalent
to 
 \be M (\Rh - q \ty\ty \II \ot \II\ty)(x_i \ot x_j)
                                                     = C^q_{i,\ty j} \, c
                       \quad\mbox{for all $i,j \in I$} \,. \label{relKul} \ee
In the purely bosonic case, this neat formula is due to Kulish \cite{Kul}. By
multiplying this equation by $\Cqi_{i,\ty j}\ty\ty$, summing over $i,j\ty$,
and using some formulae given in Ref.~\cite{Srm}, we can easily rederive
Eqn.~\reff{canima}.

In view of the conciseness of Eqn.~\reff{relKul}, we might want to use this
formula to define the algebra $\Arq$. Then we would have to show that the
Eqns.~\reff{delx} define on $\Arq$ the structure of a graded right
$\SP$--comodule algebra. This will be the case if we can show that the
relations \reff{relKul} stem from a graded subcomodule of $T(V)$, in the way
described in Section 4. In order to prove this, we consider the linear map
 \[ F : T_2(V) \lra T_0(V) \ot T_2(V) \;, \]
defined by
 \[ F(u) = \Rh(u) - q \ty\ty u - c \ty\ty \bt^q(u)
                                  \quad\mbox{for all $u \in T_2(V)$} \;. \] 
Since $\Rh$ and $\bt^q$ are homomorphisms of graded comodules, so is $F$, and
its image is the graded subcomodule of $T(V)$ we are looking for.

We close this section by a series of remarks.

\noindent
{\em Remark 6.2.} Suppose that $\KK = \CC$ (or, more generally, that $\KK$
contains the square roots of all of its elements). Then, for fixed numbers
$m$ and $n\ty\ty$, the graded right $\SP$--comodule algebras $\Arq$ with
$c \neq 0$ are all isomorphic. (This follows by simply rescaling the
generators with a suitable overall factor.) Thus, with a slight abuse of
language, we then may call this superalgebra {\em the} $\SP$--covariant
quantum Weyl superalgebra and denote it by $W_q(n \ty|\ty m)$. Needless to
say, there are really innumerable papers on quantum oscillator (super)algebras
of all types, and well--known methods to relate the various versions (see
Ref.~\cite{CDe} for an overview). Thus the main contribution of the present
work to this field is to present an $\SP$--covariant quantum oscillator
superalgebra and to provide the mathematical basis of this concept (including
the construction of the quantum supergroup $\SP$ itself).
\vspace{1.0ex}

\noindent
{\em Remark 6.3.} The defining relations of the algebra $\Arq$ can be
simplified in various ways. For example, choose elements $r_i \in \KK\:$;
$i \in I$, such that
 \[ r_i \tty\ty r_{-i} = \Cqi_{-i,\ty i}
                                 \quad\mbox{for $-r \leq i \leq -1$} \;, \]
and define
 \[ x'_i = r_i \ty\ty x_i \quad\mbox{for all $i \in I$} \,. \]
Multiplying the $i\ty$th relation of the system (I) by
$\Cqi_{-i,\ty i}\ty\ty$, we obtain the following equivalent system for the
generators $x'_i\,$:
 \[ q^{-1} \ty x'_i \tty\ty x'_{-i}
                       - \si_i \ty\ty q^{\si_i} x'_{-i} \tty\ty x'_i
        \,=\, (q - q^{-1}) \sum_{j < i} x'_{-j} \ty\ty x'_j - q^{2d} \ty c
                                \quad\mbox{for $-r \leq i \leq -1$} \;. \]
Obviously, for $i \neq -j$ the generators $x'_i$ and $x'_j$ satisfy the same
relation as $x_i$ and $x_j\ty\ty$. On the other hand, the system (II) is more
complicated for the generators $x'_i\ty\ty$. Of course, we could also
simplify the system (II), at the price that now the system (I) becomes more
complicated.

\noindent
{\em Remark 6.4.} In the case $n = m = 1\ty$, a $\Uqtt$--covariant quantum
Weyl superalgebra of the type considered in this section has been constructed
in Ref.~\cite{TWa}. However, contrary to the claim of the authors, they
consider the other version of $\Uqtt$, namely, the one associated to a basis
of the root system consisting of two odd roots. For this version of $\Uqtt$,
I have gone through all the steps of Ref.~\cite{Srm} and of the present work,
and finally have found (after adequate adjustments) the quantum Weyl
superalgebra of Ref.~\cite{TWa}.

\noindent
{\em Remark 6.5.} It should be clear that results similar to those obtained
above can also be derived in the case of graded left $\SP$--comodule algebras
(and I have done that almost completely). Actually, most of the necessary
technical tools have been mentioned in the present paper. Thus the interested
reader should have no difficulties to carry out the details. He/she might
then proceed to study the duality between the two classes of comodule
superalgebras thus obtained.

\sect{Discussion \vspace{-1ex}}
In the following, it will be useful to include the undeformed case $q = 1$ in
our discussion. The corresponding Hopf superalgebra $\SPu$ is easily
described. For $q = 1\ty$, we have $R = \id_{V \ot V}\ty\ty$, hence $A(R)$
is the universal supercommutative algebra, generated by elements $t_{ij}\tty$;
$i,j \in I$, which are homogeneous of degree $\eta_j - \eta_i$\ty\ty. The
Hopf superalgebra $\SPu$ is defined by requiring in addition the relations
\reff{invbSP} and \reff{invaSP} (with $q = 1$), or their matrix equivalents
(see Remark 5.2).

A look at the supergroup $\SPu$ helps to understand how the quantum
supergroup $\SP$ should be interpreted. The Lie group ``contained'' in
$\SPu$ is equal to ${\rm SP}(2n) \ti {\rm O}(2m)$. Correspondingly, the
algebra of functions $\SPu$ contains an element $B = B_1\ty$, the
Berezinian = superdeterminant, such that $B^2 = 1\ty$, but $B \neq 1\ty$.
It is important to recall that the Berezinian is multiplicative. By setting
$B$ equal to one, we obtain the algebra of functions on a supergroup, whose
associated Lie group is equal to ${\rm SP}(2n) \ti {\rm SO}(2m)$, and which,
accordingly, might be called $\SPSOu$.

It is tempting to conjecture that something analogous should be true in the
deformed case. More precisely, there should exist an element
\vspace{0.2ex}
$B_q \in \SP$, which is group--like, homogeneous of degree zero, central, and
satisfies $B_q^2 = 1$ but $B_q \neq 1\ty$. (We might even hope to find a
preimage of $B_q$ in $A(R)$ with analogous properties.) The quotient of $\SP$
modulo the graded Hopf ideal generated by $B_q - 1$ would then be a Hopf
superalgebra $\SPSO$. Moreover, we expect that the Hopf superalgebra pairing
$\vp$ given by Theorem 3 satisfies
 \[ \vp(B_q - 1,h) = 0 \quad\mbox{for all $h \in \Uq$} \;, \]
which implies that $\vp$ induces a Hopf superalgebra pairing $\vp_0$ of
$\SPSO$ and $\Uq$. In particular, the pairing of Theorem 3 would not satisfy
\reff{chiinj}, and the best we could hope for is that $\vp_0$ would satisfy
the implication analogous to \reff{chiinj}. Similarly, we also expect that
the universal $r$--form $\vr$ on $\SP$ described in Theorem 2 induces a
universal $r$--form on $\SPSO$.

The foregoing is closely related to the unsolved problem of whether, in the
case $n = m = 1\ty$, the $\Uqtt$--submodules $V_4$ and $\Vf$ are 
$\mbox{SPO}_q(2 \ty|\ty 2)$--subcomodules of $V \ot V$. It is easy to see
that $V_4$ (resp.~$\Vf$) is an $\mbox{SPO}_q(2 \ty|\ty 2)$--subcomodule of
$V \ot V$ if and only if $t^2_{2,-2} = 0$ (resp.~$t^2_{-2,\ty 2} = 0$). (The
fact that $V_4$ and $\Vf$ are $\Uqtt$--submodules of $V \ot V$ implies that
the matrix elements $\pi_{2,-2} = \chi(t_{2,-2})$ and 
$\pi_{-2,\ty 2} = \chi(t_{-2,\ty 2})$ satisfy the equation 
$\pi^2_{2,-2} = \pi^2_{-2,\ty 2} = 0$ in $\Uqtt^{\circ}$.) It can be shown
that in $\mbox{SPO}_q(2 \ty|\ty 2)$ (even in $A(R)$) we have
 \[ t^2_{2,-2} \ty\ty t^2_{-2,-2}
        \,=\, t^2_{2,-2} \ty\ty t^2_{2,\ty 2}
              \,=\, t^2_{-2,-2} \ty\ty t^2_{-2,\ty 2}
                    \,=\, t^2_{2,\ty 2} \ty\ty t^2_{-2,\ty 2}
    \,=\, 0 \;. \]
In the case $q = 1\ty$, it is the relation $B_1 = 1$ which implies that
$t_{-2,-2}$ and $t_{2,\ty 2}$ are invertible and hence that
$t^2_{2,-2} = t^2_{-2,\ty 2} = 0\ty\ty$. Accordingly, we conjecture that
$V_4$ and $\Vf$ are $\mbox{SPSO}_q(2 \ty|\ty 2)$--subcomodules but not
$\mbox{SPO}_q(2 \ty|\ty 2)$--subcomodules of $V \ot V$.

The construction of the quantum Berezinian $B_q$ (if it exists) is a difficult
problem. One may hope to solve it by use of a suitable Koszul complex. In the
undeformed case, the method has been described in Ref.~\cite{MKo}, the
deformed case has been discussed in Ref.~\cite{LSu}. Note, however, that in the
latter reference, the case of Iwahori, Hecke type $R$--matrices is considered,
whereas here we are in the Birman, Wenzl, Murakami setting.
 
We close this discussion by some general remarks. In the present work, we
have studied the symplecto--orthogonal quantum supergroups $\SP$ and their
graded comodule algebras. Actually, the results of the Sections 3 and 4 have
been derived under much more general assumptions. Moreover, also the
discussion of the Sections 5 and 6 should be applicable to more general cases.
A look into Section 5 shows that only some general properties of the
$R$--matrix have been used. Thus, once the $R$--matrices of (some version of)
the orthosymplectic quantum superalgebras in the vector representation have
been calculated, it should be easy to extend the results of the Sections 5
and 6 to these cases.

The reader may have noticed that I have been careful not to use the duality
of $\SP$ and $\Uq$ beyond the extent to which it is really established (see
Theorem 3). This should help to investigate this duality more carefully.

We have not given a basis of the quantum Weyl superalgebra $\Arq$. Using the
diamond lemma (see Ref.~\cite{Ber}), it should not be difficult to prove that
(with respect to some suitable ordering of the indices) the ordered monomials
in the generators $x_i$ form a basis of this algebra; of course, the exponent
of $x_i$ should be restricted to $\{0\tty,1\}$ if $x_i$ is odd.
\vspace{3.0ex}

\noindent
{\bf Acknowledgements} \\
Throughout the investigations which led to the present work, I have had
various highly fruitful discussions with Ruibin Zhang on most of the topics
considered here. I am very grateful to him for sharing with me his insights
into the theory of quantum supergroups. Part of this research has been
carried out during a visit of the author to the Department of Mathematics at
the University of Queensland. The kind invitation by Mark Gould and the
hospitality extended to the author by the members of the department are
gratefully acknowledged.

\end{document}